\documentclass[noams]{compositio}

\usepackage{amsbsy,amssymb,amscd,amsfonts,latexsym,amstext,delarray,amsmath,stackrel,lineno,tabularx,url,MnSymbol,tikz,cite,yfonts}
\usepackage{float} 

\usepackage[pdfborder={0 0 0}]{hyperref}
\usepackage{mathrsfs}
\hypersetup{
    colorlinks = true,
    linkcolor = blue,
    anchorcolor = red,
    citecolor = blue,
    filecolor = red,
    pagecolor = red,
    urlcolor = magenta
}

\usepackage[final]{showlabels}

\showlabels{cite}

\usepackage{graphicx}

\usepackage{pdfpages}

\usepackage{blindtext}
\usepackage{etoolbox}
\newcommand{\printthis}[2][true]{%
\ifbool{#1}{%
#2%
}{%
}%
}

\newtheorem{thm}{Theorem}[section] 
\newtheorem{defn}[thm]{Definition} 

\newtheorem{prop}[thm]{Proposition}
\newtheorem{cor}[thm]{Corollary}
\newtheorem{lem}[thm]{Lemma}

\newtheorem{rem}[thm]{Remark}

\def\fii{{\rm  finite}}

\def\Hom{{\rm Hom}}
\def\id{{\rm id}}

\def\Spec{{\rm Spec\,}}

\def\B{{\mathbb B}}

\def\F{{\mathbb F}}

\def\N{{\mathbb N}}

\def\Q{{\mathbb Q}}
\def\R{{\mathbb R}}
\def\Z{{\mathbb Z}}

\def\pii{\pi^{(2)}}
\def\pij{\pi^{\rm comb}}
\def\pik{\pi^{\rm new}}

\def\cC{{\mathcal C}}

\def\cJ{{\mathcal J}}

\def\cM{{\mathcal M}}
\def\cO{{\mathcal O}}

\def\cS{{\mathcal S}}

\def\qqq{\,,\,~\forall}

\newcommand{\ie}{{\it i.e.\/}\ }
\newcommand{\eg}{{\it e.g.\/}\ }

\newcommand{\opcit}{{\it op.cit.\/}\ }

\def\dellt{{\Delta'}}

\def\sin{{{\rm sin}}}

\def\id{{\mbox{Id}}}

\def\ker{{\mbox{Ker}}}

\def\Hom {{\mbox{Hom}}}

\def\sss{{\mathfrak s}}
\def\gop{{\Gamma^{\rm op}}}

\newcommand*{\xhat}[1]{#1\kern-0.08em\widehat{\phantom{#1}}}

\def\Se{\frak{Sets}}
\def\set2{\frak{Sets}^{(2)}}
\def\fin{\frak{Fin}}

\def\vvert{{\Vert}}
\def\hq{{{\cO\subset H\Q}}}

\def\dop{{\Delta^{\rm op}}}

\def\Ses{{\Se_*}}
\def\sses{{\mathcal S_*}}
\def\gsses{{\Gamma\cS_*}}
\def\gse{\underline{\Gamma\Se}_*}

\def\gamm{\widehat{\Gamma}}

\def\smod{{\sss-{\rm Mod}}}
\def\gqu{{$\Gamma$-2-set }}
\def\gqus{{$\Gamma$-2-sets }}

\title{$\overline{\Spec\Z}$ and the Gromov norm}
\author{Alain Connes}
\email{alain@connes.org}
\address{I.H.E.S. and Ohio State University}
\author{Caterina Consani}
\email{kc@math.jhu.edu}
\address{Department of Mathematics, The Johns Hopkins
University\newline Baltimore, MD 21218 USA}

\thanks{The second author is partially supported by the Simons Foundation collaboration grant n. 353677.}
\begin{document}
\begin{abstract}
We define the homology of a simplicial set with coefficients in a Segal's $\Gamma$-set ($\sss$-module). We show the relevance of this new homology with values in $\sss$-modules by proving that taking as coefficients the $\sss$-modules at the archimedean place over the structure sheaf on $\overline{\Spec\Z}$ as in \cite{CCprel}, one obtains on  the singular homology with real coefficients of a topological space $X$, a norm equivalent to the Gromov norm. Moreover, we prove that the two norms agree when $X$ is an oriented compact Riemann surface.
\end{abstract}
\maketitle

\vspace*{6pt}\tableofcontents 

\section{Introduction}
The notion of a $\Gamma$-set is a fundamental constituent in mathematics: it is the most embracing generalization of the datum given on a set by a commutative addition with a zero element and it  provides a common framework for many of the present efforts to understand the ``field with one element". In \cite{CCprel} we defined on the Arakelov compactification $\overline{\Spec\Z}$ of the algebraic spectrum of the integers a structure sheaf of $\Gamma$-rings which agrees with the classical structure sheaf when restricted to $\Spec \Z$, but whose stalk at the archimedean place uses in a crucial way the new freedom of moving from the  category of abelian groups to that  of $\Gamma$-sets. To define $\Gamma$-sets one first introduces  the small, full subcategory $\Gamma^{\rm op}$  of  the category $\fin_*$ of pointed finite sets, whose objects are pointed sets $k_+:=\{0,\ldots ,k\}$, for each integer $k\geq 0$ ($0$ is the base point) and with morphism the sets $\gop(k_+,m_+)=\{f: \{0,1,\ldots,k\}\to\{0,1,\ldots,m\}\mid f(0)=0\}$. A $\Gamma$-set is then defined as a (covariant) functor $\gop\longrightarrow\Ses$ between pointed categories and 
the morphisms  in this category  are natural transformations. The closed structure of the category $\Gamma\Ses$  of $\Gamma$-sets  is defined by setting
\begin{equation}\label{closedstructure0}
	\gse(M,N)=\{k_+\mapsto \Gamma\Ses(M,N(k_+\wedge -))\},
\end{equation}
where $\wedge$ is the smash product of pointed sets. This formula uniquely defines the smash product of $\Gamma$-sets by applying the adjunction 
$$
\gse(M_1\wedge M_2,N)=\gse(M_1,\gse(M_2,N)).
$$
The notions of rings and modules then acquire a meaning in this symmetric monoidal closed category. In particular, 
 $\Gamma$-sets can equivalently be viewed as modules over the simplest $\Gamma$-ring $\sss: \gop\longrightarrow\Ses$ whose underlying $\Gamma$-set is the  identity functor, whence the name $\sss$-module to denote a $\Gamma$-set,  and the more suggestive notation for morphisms in $\Gamma\Ses$
 $$
 \Hom_\sss(M,N):=\Gamma\Ses(M,N), \ \ \underline\Hom_\sss(M,N):= \gse(M,N).
 $$
 Abelian groups form a full subcategory of the category $\smod$ of $\sss$-modules: the inclusion functor associates to an abelian group $A$ the functor (Eilenberg-Mac Lane object) $HA:\gop\longrightarrow\Ses$ which  assigns to a finite pointed set $X$ the pointed set of $A$-valued maps  on $X$ vanishing at the base point of $X$ (\!\!\cite{DGM}, 2.1.2). \newline 
  At the conceptual level, it is important to make as explicit as possible the link between the category $\smod$ and the naive interpretation of  vector spaces over $\F_1$ as pointed sets (see \cite{KS}). This link can be understood by viewing $\sss$-modules as {\em pointed objects} in the topos  $\gamm$ of covariant functors  $\gop\longrightarrow\Se$. Thus,  provided one works in $\gamm$, one may think of our basic objects as ``pointed sets". The reason for this choice of topos is to provide room for the identity functor $\id: \gop\longrightarrow \gop$ which defines   the simplest $\Gamma$-ring: $\sss$. In other words both $\gop$ and objects in $\gamm$ are based on the idea of pointed sets which underlies the naive interpretation of $\F_1$. In this way one reaches a workable framework that extends strictly the category of $\Z$-modules.
  
 To perform homological algebra one needs, guided by the Dold-Kan correspondence, to move from the basic category $\smod$ to its simplicial version, namely the category $\gsses$ of $\Gamma$-spaces, where $\sses$  denote the category of pointed simplicial sets, \ie contravariant functors $\Delta \longrightarrow \Se_* $ where $\Delta$ is the ordinal number category and $\Se_*$ is the category of pointed sets. The category $\gsses$  plays a  central role in 
 \cite{DGM}.  We denote by  $\underline\Hom_\sses$  the internal hom functor in $\sses$. As explained in \opcit, one can use the closed structure of $\sses$  to endow   $\gsses$  with the structure of a symmetric monoidal closed category.  The closed structure is defined as follows 
 \begin{equation}\label{closedstructure}
	\underline\Hom_\gsses(M,N):=\{(k_+,[q])\mapsto \Hom_\gsses(M\wedge \Delta[q]_+,N(k_+\wedge -))\}.
\end{equation}
The monoidal structure is given by the smash product where $M\wedge N$ is defined using the  closed structure and can be described as a Day's product (see \opcit 2.1.2.1)
 $$
 (M\wedge N)(Z)=\int^{(X,Y)}\left( M(X)\wedge N(Y)\right)\wedge \gop(X\wedge Y,Z).
 $$ 
 The key result of Lydakis (\opcit Theorem 2.1.2.4) states that there are choices of coherency isomorphisms so that the triple $(\gsses, \wedge,\sss)$ is a symmetric monoidal closed category.
 
Our goal  is to use $\Gamma$-spaces to perform homological algebra in the category $\smod$ by applying an analogue of the Dold-Kan correspondence. 
 For our arithmetic applications it is crucial to work with {\em non-fibrant} $\Gamma$-spaces  and  define a suitable substitute for the homotopy groups. In homotopy theory the Kan extension property is used in two ways:\newline
 -~to show that the relation of homotopy is an equivalence relation,\newline
 -~to define the group structure on $\pi_n$ for $n\geq 1$.	\newline
To define the homology  $H_n(X,F)$ of a pointed simplicial set $X$ with coefficients in an $\sss$-module $F$, the problem to obtain the substitute of the group structure does not arise since, already in the classical case where $F=HA$ corresponds to an abelian group, the interchange law shows that the group structure in homology is the same as that inherited from the underlying $\Gamma$-set (see Remark \ref{homotopgroup}).  Thus the issue created by the lack of the Kan extension property occurs mainly at the level of pointed (non fibrant) simplicial sets $X$, and $\gop$ is not involved there.  One thus needs, as an intermediate step,  to extend the combinatorial construction of the homotopy $\pi_n(X,\star)$ for a pointed simplicial set which is not fibrant. This step is described in 
  Section \ref{sectcombinatorial1} of the present paper.  The main difficulty to obtain a meaningful combinatorial notion  is that the relation of homotopy between $n$-simplices $x,y\in X_n$ as in \cite{May} Definition 3.1 is no longer an equivalence relation. By definition (see \opcit)
 \begin{equation}\label{Maydefn}
  R=\{(x,y)\in X_n\times X_n\mid \partial_j x=\partial_jy\, \forall j \, \&\,  \exists z \mid \partial_jz=s_{n-1}\partial_j x\, \forall j<n, \ \partial_nz=x,\, \partial_{n+1}z=y\}
  \end{equation}
  The simplices involved in the definition of $\pi_n$  correspond  to the elements of  $\Hom_\sses(S^n,X)$, \ie by Yoneda's lemma to $x\in X_n$ with $ \partial_j x=*\, \forall j$. Here   $S^n$ is the combinatorial sphere, \ie  the pointed simplicial set $(\Delta[n],\partial \Delta[n])$ obtained by collapsing the boundary $\partial \Delta[n]$ of the standard simplex to a single base point. The  relation $R$ on $\Hom_\sses(S^n,X)\subset X_n$ coincides with the relation on the $0$-skeleton $Y_0$ associated to the two boundary maps $\partial_j:Y_1\to Y_0$, where $Y:=\Omega^n(X)$ is obtained from $X$ by iterating $n$-times the endofunctor $\Omega:\sses\longrightarrow \sses$ of \cite{Moore} (Definition 1.6).  In this way one reduces the problem to the definition of $\pi_0Y$ for $Y=\Omega^n(X)$. Then one can simply define $\pi_0Y$ as the quotient $\pij_0Y$ of $Y_0$ by the equivalence relation generated by $R$.  
  \begin{defn}\label{defnpij} Let $n\geq 0$ be an integer and  $X$ a  pointed simplicial set.  Define
 \begin{equation}\label{pins00}
\pij_n(X):=\pij_0(\Omega^n(X))=\Hom_\sses(S^n,X)/\tilde R
\end{equation}
where $\tilde R$ is the equivalence relation generated by the restriction of the relation $R$ of \eqref{Maydefn}.
 \end{defn}
 This notion developed in Section \ref{sectcombinatorial1} suffices for the goals of the present paper,  but for future applications we also wish to keep the finer information contained in the relation $R$. 
  This is achieved by introducing the topos $\set2$ in which the finer notion, denoted  $\pii_n(X,\star)$,  takes its value : \ie  $\pii_n(X,\star)$ is a $2$-set, \ie  an object of $\set2$. This construction is described in Section \ref{sectcombinatorial1bis} where we also show that the topos $\set2$ is related to the topos of quivers.   In Section \ref{secthomotop} we then obtain  a	 general definition of homology of $\Gamma$-spaces, considered as simplicial $\Gamma$-sets. This homology is not, in general,   a group but is a  \gqu \ie  a pointed covariant functor $\gop \longrightarrow \set2_*$.

 In Section \ref{secthomology} we construct, given an integer $n\geq 0$, an arbitrary pointed simplicial set $X$ and an arbitrary $\sss$-module ($\Gamma$-set) $F$,  the homology $H_n(X,F)$ as follows
  \begin{defn}\label{defnhomol} Let $n\geq 0$ be an integer,  $X$ a  pointed simplicial set, and $F$ an $\sss$-module.   Define 
 \begin{equation}\label{pins1}
H_n(X,F):=\{k\mapsto \pij_n(F\circ ( X\wedge k_+ ))\}
\end{equation}
as an $\sss$-module.
 \end{defn} 
  As in \cite{DGM}, we extend the $\Gamma$-set $F$ to an endofunctor of the category of pointed sets.
 When $F=HA$ for an abelian group $A$, $H_n(X,F)$  coincides with the standard definition of homology :
 \begin{thm}\label{propcompareintro}
	Let $A$ be an abelian group, and $X$ a pointed simplicial set. For any integer $n\geq 0$ one has the equality of $\sss$-modules
	\begin{equation}\label{equbasic}
  H_n(X,HA)=H(H_n(X,A))
\end{equation} 
where $H_n(X,A)$ is the (reduced) abelian group homology  of $X$ with coefficients in $A$.
\end{thm}
   Again, we stress the fact that we apply Definition \ref{defnhomol} in cases where the pointed simplicial set $F\circ ( X\wedge k_+ )$ is {\em not fibrant}. In particular, in our applications the $\sss$-modules $H_n(X,F)$ are rarely groups. 
 
 In Section \ref{sectgromov} we apply Definition \ref{defnhomol} to the $\sss$-modules  we introduced in  \cite{CCprel}, at the archimedean place of $\overline{\Spec\Z}$. We show that these coefficients yield a semi-norm on the ordinary singular homology $H_n(X,\R)$ of a topological space $X$ and our goal is to compare this semi-norm with the Gromov norm, whose definition is recalled in Section \ref{simplicialvol}. In Section \ref{sectspecZ} we review our construction 
 (see \cite{CCprel}) of the structure sheaf $\hq$ of $\sss$-algebras on $\overline{\Spec\Z}$.  The  sheaves $\cO(D)$ associated to Arakelov divisors $D=D_\fii+D_\infty$, as  in \opcit   provide a one parameter family of $\sss$-modules $\Vert H\R \Vert_\lambda$ ($\lambda\in\R_+$) which we can use as coefficients in formula \eqref{pins1}. In Section \ref{sectequivalence}, Proposition \ref{conj}, we prove that for any  topological space $X$ the filtration  of the singular homology group $H_n(X,\R)$ by the $H_n(X,\Vert H\R \Vert_\lambda)$ defines a semi-norm 
which is equivalent to the Gromov norm. 

The final Section \ref{sectequal} is entirely devoted to show that the two norms on $H_n(X,\R)$: the Gromov norm and our new norm,  are in fact equal when $X=\Sigma$ is a compact Riemann surface. The difficulty in the proof of this result is due to  the fact that in order to obtain elements of the homology $H_2(\Sigma,\Vert H\R \Vert_\lambda)$ one needs to get singular chains which are not only cycles but are such that all their simplicial boundaries actually vanish. While one knows that this Moore normalization is possible the problem is to effect it without increasing the $\ell^1$-norm of the chain : this requires a delicate  geometric work described in Sections \ref{block} and \ref{block1}. One then obtains the desired equality in the form of the following
   \begin{thm}\label{chainthmintro} Let $\Sigma$ be a compact Riemann surface and $[\Sigma]$ its fundamental class in homology. Then $[\Sigma]$ belongs to the range of the canonical map $H_2(\Sigma,\Vert H\R \Vert_\lambda)\to H_2(\Sigma,\R)$ if and only if $\lambda$ is larger than the Gromov norm of $[\Sigma]$. 	 	
 	 \end{thm}
We expect that a similar statement holds in hyperbolic geometry in any dimension. The natural test ground for the homology 
$H_n(X,\Vert H\R \Vert_\lambda)$ is in  hyperbolic spaces since the Gromov norm does not vanish there for $n>1$  while it vanishes identically on all spheres. This is in contrast with the construction of the spectra associated to $\Gamma$-spaces $M$ where the associated  endofunctor $X \mapsto M\circ X$ is only tested on spheres.

 \section{Homology of a simplicial set with coefficients in an $\sss$-module}\label{sectcombinatorial}
 
  Our goal in this section is to reach a good definition of the homology of a pointed simplicial set with coefficients in an $\sss$-module and  to show that it generalizes the standard notion in algebraic topology.  This is achieved in Definition \ref{fact1} and Theorem \ref{propcompare}.  As a preliminary step we need to refine the definition of the homotopy groups $\pi_n$ by remaining at the combinatorial level and ignoring the group structure. Classically (see \eg\cite{DGM} Appendix A.2.3), the function space of maps  between pointed simplicial sets $X$ and $Y$ is defined as the pointed simplicial set:
\begin{equation}\label{functionspace}
	 {\rm Map}_*(X,Y):=\underline\Hom_\sses(X,\sin\vert Y\vert)
\end{equation}
This amounts to replace $Y$ with the {\em fibrant} simplicial set $\sin\vert Y\vert$ and it entails that the $\pi_n$,  defined using such a fibrant replacement (see also \opcit A.2.5.1),  are then groups for $n\geq 1$ (abelian for $n>1$). Thus in the definition of the homotopy groups of a  $\Gamma$-space $M:\gop\longrightarrow\sses$ (see \opcit Definition 2.2.1.2)
 \begin{equation}\label{functionspace1}
	 \pi_qM:=\varinjlim_k \pi_{k+q}M(S^k)
\end{equation}
 the terms involved in the colimit are groups, hence $\pi_qM$ is an abelian group. \newline
For our applications however, the simplification effected by the definition \eqref{functionspace}  hides certain finer features of $\Gamma$-spaces which become relevant for  arithmetic constructions.  We shall thus work directly in the category $\gsses$ without performing this fibrant replacement.
 
 \subsection{Homotopy  for pointed simplicial sets}\label{sectcombinatorial1}
 In order to define the new homotopy $\pik_n(X,\star)$ for a general pointed simplicial set $(X,\star)$, we shall first reduce to the case of $\pik_0$, following \cite{Moore} Definition 1.9. One defines an endofunctor $\Omega$ of $\sses$ which associates to a pointed simplicial set $(X,*)$ the pointed simplicial set $\Omega(X,*)$ defined as follows (with $k$ a positive integer)
 \begin{equation}\label{omegadef}
\Omega(X,*)_k:=\{x\in X_{k+1}\mid \partial_0(x)=*, \ \partial_{i_0}\ldots \partial_{i_k}x=*\qqq i_j\in \{0,\ldots ,k+1\}\}
\end{equation}
with the simplicial structure given by faces
\begin{equation}\label{omegadef1}
\partial_j:\Omega(X,*)_k\to \Omega(X,*)_{k-1}, \  \ \partial_j(x)=\partial_{j+1}^X(x)
\end{equation}
and degeneracies
\begin{equation}\label{omegadef2}
s_j:\Omega(X,*)_k\to \Omega(X,*)_{k+1}, \ \  s_j(x)=s_{j+1}^X(x).
\end{equation}
  The definition of the homotopy $\pik_n(X,\star)$  is then reduced to that of $\pik_0$ for the simplicial set $\Omega^n(X)$ obtained after iterating the endofunctor $\Omega$ n-times : 
\begin{equation}\label{pindef}
\pik_n(X,\star):=\pik_0(\Omega^n(X)).
\end{equation}
One shows by induction\footnote{Note that a product $\partial_{i_0}\ldots \partial_{i_k}$ can be reordered using the simplicial rules so that the indices fulfill $i_0\geq i_1\geq \ldots \geq i_k$.} on $n$ that 
  \begin{equation}\label{omegadef3}
\Omega^n(X,*)_k =\{x\in X_{n+k}\mid \partial_j(x)=*, \ \forall j<n, \ \partial_{i_0}\ldots \partial_{i_k}x=*\qqq i_j\in \{0,\ldots ,k+n\}\}
\end{equation}
while the face and degeneracies are obtained as in \eqref{omegadef1} and \eqref{omegadef2} but using $\partial_{j+n}^X$ and  $s_{j+n}^X$. 
One  describes directly the first levels of $\Omega^n(X)$ as follows
\begin{lem}\label{levels} Let $(X,*)$ be a  pointed simplicial set.\newline
$(i)$~The $0$-skeleton $(\Omega^n(X))_0$ is the set of simplices $x\in X_n$  with all  $\partial_j(x)$  equal to the base point.\newline
$(ii)$~$(\Omega^n(X))_0$ coincides with ${\Hom}_\sses(S^n,X)\subset X_n$ where   $S^n$ is obtained by collapsing the boundary $\partial \Delta[n]$ of the standard simplex to a single base point\footnote{This is not the definition used in \cite{DGM}, where  $S^n$ is defined as the $n$-fold smash product $S^1\wedge \dots \wedge S^1 $ of $S^1=\Delta[1]/ \partial \Delta[1]$. This distinction in the definition of the homotopy groups is irrelevant in the fibrant case since the geometric realizations are homeomorphic, but as in \cite{Moore} our choice is more convenient to compute the set of maps using Yoneda's Lemma.}.\newline
$(iii)$~The $1$-skeleton $(\Omega^n(X))_1$ is the set of $x\in X_{n+1}$ which fulfill the conditions
$$
\partial_i\partial_j(x)=*\qqq i,j, \ \ \partial_j(x)=*\qqq j\in \{0,\ldots, n-1\}.
$$
$(iv)$~The boundaries $\partial_i:(\Omega^n(X))_1\to (\Omega^n(X))_0$ for $i=0,1$ are given by $\partial_n$ and $\partial_{n+1}$.\newline
$(v)$~The relation $R$ on $(\Omega^n(X))_0={\Hom}_\sses(S^n,X)\subset X_n$ given by 
 $$
x R y \iff \exists z\in (\Omega^n(X))_1~\text{s.t.}~ \partial_0 z=x~\text{and}~ \partial_1 z=y
 $$ 
coincides with the relation of homotopy between $n$-simplices  as in \eqref{Maydefn}. 
 \end{lem}
 \proof $(i)$~Follows from \eqref{omegadef3} for $k=0$.\newline
$(ii)$~By Yoneda's lemma one checks that the morphisms $y\in{\Hom}_\sses(S^n,X)$ \ie the elements of ${\Hom}_\sses(\Delta[n],X)$ which send $\partial \Delta[n]$ to the base point, are the same as the elements of the $0$-skeleton $(\Omega^n(X))_0$.\newline
$(iii)$~Follows from \eqref{omegadef3} for $k=1$.\newline
$(iv)$~Follows from $\partial_j=\partial_{j+n}^X$ for $j=0,1$.  \newline
$(v)$~This follows from the previous part of the lemma since the relation \eqref{Maydefn} restricts to 
 \begin{equation}\label{Maydefnbis}
  R=\{(x,y)\in X_n\times X_n\mid \partial_j x=\partial_jy=*\, \forall j \, \&\,  \exists z \mid \partial_jz=*\, \forall j<n, \ \partial_nz=x,\, \partial_{n+1}z=y\}
  \end{equation}
 \endproof
 \begin{rem}\label{geomrem} The geometric meaning of the endofunctor $\Omega$ can be understood as explained to us by B. Dundas: first, as in \cite{DGM} A.2.7, one has a combinatorial model $PX$ mimicking the path space of a simplicial set $X$ by precomposing the functor $X$ with the endofunctor $[0]\coprod \bullet: \Delta\longrightarrow \Delta$. This simply shifts the indices \ie one has $(PX)_k=X_{k+1}$ and the indices of faces and degeneracies are shifted by $1$. The link with ordinary paths is given by precomposing with the  morphism of simplicial sets
 $
\gamma: \Delta[1]\times \Delta[q]\to \Delta[q+1]
 $, 
 associated as $\gamma:=N(p)$ by the nerve functor $N$ to 
 $$
 p: [1]\times [q]\to [q+1], \ \ p(0,j):=0 \ \forall j, \ \ p(1,j):=j+1 \ \forall j
 $$
 Requiring that the two end points of the path associated to  $x\in X_{k+1}=\Hom_\Delta(\Delta[q+1],x)$ are equal to the base point $*$ (when $X$ is pointed) gives exactly the conditions of \eqref{omegadef} defining $\Omega(X)$.	When $X$ is fibrant one obtains in this way a model for its loop space.
 \end{rem}

      For a {\em fibrant} simplicial pointed set $X$, the relation \eqref{Maydefnbis} is an {\em equivalence} relation and the quotient by this relation defines $\pi_0(\Omega^n(X))$ which is known to be a group, for $n\geq 1$ (see \cite{May,Moore}, or  Theorem 7.2 in Chapter III of \cite{GJ}).
    Note also that when $X$ is fibrant the above equivalence relation on ${\Hom}_\sses(S^n,X)\subset X_n$ coincides with the one defined by the two boundary maps from the $1$-skeleton of the simplicial set $\underline\Hom_\sses(S^n,X)$ (see \cite{Moore} Lemma 1B.3).    
    \newline
On the other hand, the simplicial sets $X$ we consider here are {\em not necessarily fibrant}  and the relation $R$  is not  in general transitive (nor symmetric). The easy solution to bypass this problem is to define $\pij_n(X)$ as the quotient  by the equivalence relation generated by the relation $R$ in agreement with Definition \ref{defnpij}. This provides a first notion of homotopy which  suffices for the goal of the present paper. One has by construction
 \begin{equation}\label{pij}
\pij_n(X,\star):=\pij_0(\Omega^n(X)).
\end{equation}	
We state simple properties of this combinatorial notion
\begin{prop}\label{propprepa} $(i)$~Let $X$ be a pointed simplicial set and $k>0$ an integer. Then for any $n$ 
$$
\pij_n(X \wedge k_+)=\pij_n(X)\wedge k_+.
$$
$(ii)$~Let $X,Y$ be pointed simplicial sets, one has for any $n$
$$
\pij_n(X \times Y)=\pij_n(X)\times \pij_n(Y).
$$	
\end{prop}
\proof $(i)$~An element $x\in (X \wedge k_+)_n$, $x\neq *$ is of the form $x=(a,j)$ with $a\in X_n$ and $0<j\leq k$. Two elements $x=(a,j)$ and $x'=(a',j')$ fulfill $(x,x')\in R$ as in  
\eqref{Maydefnbis} if and only if  $j=j'$ and  $(a,a')\in R_X$ since the boundaries preserve the index $j$. \newline
$(ii)$~This follows since $(X\times Y)_n=X_n\times Y_n$ and the boundaries act componentwise. \endproof 
\subsection{The finer notion $\pii_n(X)$ and the topos $\set2$}\label{sectcombinatorial1bis}
 For later applications to Arakelov divisors Definition \ref{defnpij} is too coarse and  
 one would like to 
\begin{itemize}
 \item keep all the information about the relation $R$ and
 \item still think of $\pik_0$ as a set.	
 \end{itemize}
The idea of ``topos" of Grothendieck \cite{SGA} comes to the rescue  providing a satisfactory  answer. We consider the topos $\set2$ of contravariant functors to the category of sets from the small category obtained by restricting the objects of $\Delta$ to $[0]$ and $[1]$ and keeping the same morphisms as in the following definition 
 \begin{defn}\label{defnpin} $(i)$~Let $X:\dop\longrightarrow \Se$ be a simplicial set. We define $\pii_0(X)$ as the object of $\set2$ which is the restriction of the functor $X$ to the full subcategory of $\Delta$ with  objects $[0],[1]$ and same morphisms as $\Delta$. \newline
 $(ii)$~Let $X$ be a pointed simplicial set, then we define 	\[\pii_n(X):=\pii_0(\Omega^n(X)).
 \]
 \end{defn}
It turns out that the topos $\set2$  can also be described as the dual of the small category with a single object whose morphisms form the  monoid $\cM$ with three elements $1,m_0,m_1$ and the multiplication table  specified by the rule $m_j x=m_j$ for all $j\in \{0,1\}$. 

\begin{prop}\label{propquivers} The topos $\set2$ is the same as  the dual of the monoid $\cM$.	
\end{prop}
\proof 
By definition  an object $F$ of the topos $\set2$ is a pair of sets $F(0),F(1)$,  with two maps $\partial_j:F(1)\to F(0)$, $j\in \{0,1\}$ and a map $s:F(0)\to F(1)$ such that $\partial_j\circ s=\id$. This implies that $s:F(0)\to F(1)$ is an injection and one can thus view $F(0)$ as a subset of $F(1)$ and consider the  two self-maps $T_j=s\circ \partial_j: F(1)\to F(1)$. They  fulfill the rule  
$$
T_i\circ T_j=T_j \qqq i,j\in \{0,1\}
$$
since $s \circ (\partial_i\circ s) \circ \partial_j=s \circ (\id) \circ \partial_j=s\circ \partial_j$. Thus one obtains an object 
in the dual $\hat \cM$ of the monoid  defined by the opposite of the above rules. Conversely given an object $X$ of $\hat \cM$, \ie a set $X$ endowed with a right action of $\cM$ one defines an object of $\set2$ by setting $F(1):=X$, $F(0):={\rm Range}( T_j)$ which does not depend on the choice of $j\in \{0,1\}$. One lets $s:F(0)\to F(1)$ be the inclusion as a subset, and $\partial_j:F(1)\to F(0)$ is given by $T_j$. One checks that $\partial_j\circ s=\id$. One obtains in this way two functors $\set2\longrightarrow \hat \cM$ and $\hat \cM\longrightarrow \set2$ which are inverse of each other. \endproof 
\begin{rem} The topos $\set2$ is closely related to the topos of quivers but is not the same.
	 By definition a quiver is a directed graph where loops and multiple arrows between two vertices are allowed, \ie  a multidigraph. It is described in a precise manner by two sets $(V,E)$ which represent the vertices and the edges of the graph 
and two maps $d_j:E\to V$, $j\in \{0,1\}$, which give the source and the target of an edge. These two maps do not fulfill any condition. One obtains an object of $\set2$ by setting $F(0):=V$ and $F(1):=V \coprod E$, the disjoint union of $V$ and $E$. One then lets $s:F(0)\to F(1)$ be the inclusion while the two maps $\partial_j:F(1)\to F(0)$, $j\in \{0,1\}$ are given by $\partial_j=(\id, d_j):V \coprod E\to V$. One has $\partial_j\circ s=\id$ by construction, and $F$ is an object of $\set2$. Conversely, given an object $F$ of $\set2$ one obtains a quiver by setting 
$$
V=F(0), \ \ E=F(1)\setminus s(F(0)), \ \ d_j=\partial_j \vert E
$$
However this second construction is not functorial. In fact the topos of quivers has two points given by the functors to the set of vertices and the functor to the set of edges. Similarly these two functors give the two points of the topos $\set2$ but in the latter case the functor to the set of edges never takes the value $\emptyset$ when the functor to the set of vertices takes a non-empty value. This shows that the topos of quivers is not the same as the topos $\set2$. 
\end{rem}

 \begin{figure}[t!]	\begin{center}
\includegraphics[scale=0.6]{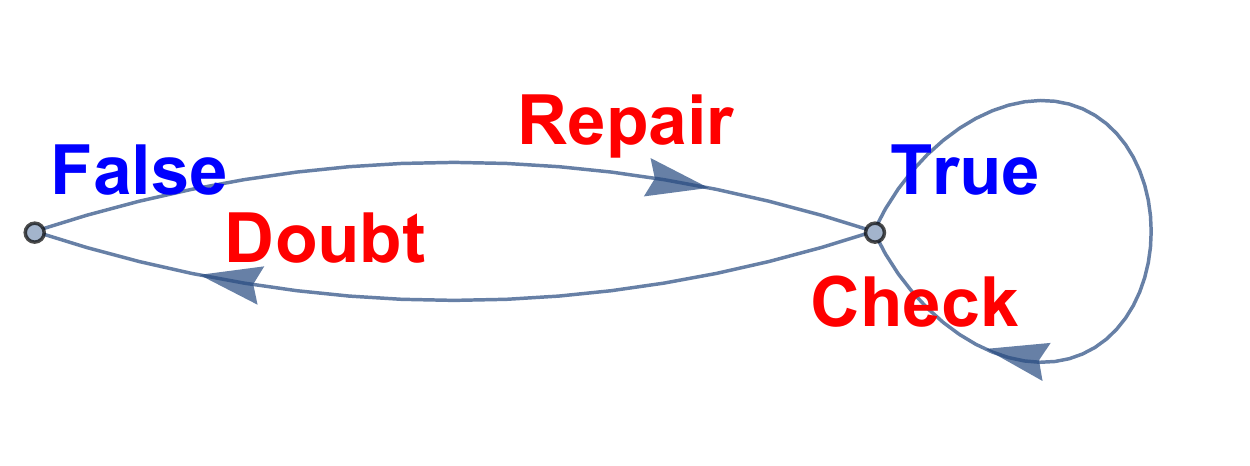}
\end{center}
\caption{Subobject classifier for $\set2$. \label{subobject} }
\end{figure} 
 \begin{lem}\label{omega} The subobject classifier  $\Omega$ of the topos $\set2$ is the object with two vertices {False, True} and five edges which besides the two degenerate ones form the  graph of Figure \ref{subobject}. 	
 \end{lem}
\proof It is a general fact (see \cite{MM}, \S I.4) that for a topos which is the dual of a monoid $\cM$ (viewed as a category with a single object), \ie the topos of sets with a right action of $\cM$, the subobject classifier is given by the set $\cJ$ of right ideals of $\cM$ on which the right action of $\cM$ is defined by 
$$
J.m:=\{n\in \cM\mid mn\in J\}\qqq J\in \cJ, \ m \in \cM.
$$
 Taking the above $\cM$ with three elements $1,m_0,m_1$ and the multiplication table specified by the rule $m_j x=m_j$ for all $j\in \{0,1\}$, one finds that $\cJ$ contains five elements
 $$
 \cJ=\{\emptyset, \{m_0\}, \{m_1\}, \{m_0,m_1\}, \cM\}
 $$
 and that the right action $T_j$ of $m_j\in\cM$  fixes $\emptyset$ and $\cM$ (which are hence degenerate edges, \ie vertices) while $T_j\{m_j\}=\cM$ and $T_i\{m_j\}=\emptyset$ for $i\neq j$. Thus the set $V$ of vertices contains two elements $\emptyset$ and $\cM$ and the non-degenerate edges are the three edges shown in Figure \ref{subobject}.\newline
 The reason for renaming the vertices $\emptyset$ as ``False" and $\cM$ as ``True" and for the choice of the labels of the edges comes from the construction of the classifying map associated to a subobject $G'$ of an object $G$ in $\set2$. One finds that the classifying map $f$ is obtained as follows as a map from $G$ to  $\Omega$:
 \begin{enumerate}
 	\item $\epsilon\in G'\Rightarrow f(\epsilon)={\rm True}$
 	\item $\epsilon\notin G' \, \text{and} \ \partial_j \epsilon\notin G'\Rightarrow f(\epsilon)={\rm False}$
 	\item $\epsilon\notin G', \,  \ \partial_0 \epsilon\notin G' \, \text{and} \ \partial_1 \epsilon\in G'\Rightarrow f(\epsilon)={\rm Repair}$
 	\item $\epsilon\notin G', \ \partial_0 \epsilon\in G' \, \text{and} \ \partial_1 \epsilon\notin G'\Rightarrow f(\epsilon)={\rm Doubt}$
 	\item $\epsilon\notin G', \ \partial_0 \epsilon\in G' \, \text{and} \ \partial_1 \epsilon\in G'\Rightarrow f(\epsilon)={\rm Check}$ .	
 \end{enumerate}
The terminology ``False" and ``True" is the standard one for the two extremes in subobject classifiers, the notations for the edges are suggestive but more arbitrary. \endproof 
This determination of the subobject classifier shows that the topos $\set2$ is two valued and not boolean (see \cite{MM}, VI).  
 
 \subsection{Homotopy of $\Gamma$-spaces}\label{secthomotop}

If $\cC$ is a pointed category with initial and final object denoted $*$, one defines (see \cite{DGM}) the category of $\Gamma$-objects of $\cC$ as the category $\Gamma\cC$ of pointed  covariant functors  $\gop \longrightarrow\cC$. This construction applies to the category $\cS_*$ of pointed simplicial sets to yield the category $\gsses$ of  $\Gamma$-spaces. It also applies to the category $\set2_*$ of pointed objects in $\set2$. We shall call \gqus the objects of $\Gamma\set2_*$.
 \begin{prop}\label{defnforgamma} $(i)$~Let $X$ be a $\Gamma$-space and $n\in \N$. Then the map 
 $
 k\mapsto \pii_n(X(k_+))
 $ (resp. $k\mapsto \pij_n(X(k_+))$) extends to a pointed  covariant functor $\pii_n(X):\gop \longrightarrow\set2_*$ (resp. to a $\Gamma$-set).\newline
 $(ii)$~For $n\in \N$, $\pii_n$ defines a functor $\pii_n:\gsses\longrightarrow\Gamma\set2_*$ from $\Gamma$-spaces to \gqus\!.\newline
 $(iii)$~For $n\in \N$, $\pij_n$ defines a functor $\pij_n:\gsses\longrightarrow\Gamma\Se_*$ from $\Gamma$-spaces to $\smod$.
 \end{prop}
\proof This follows from
 the naturality of  Definitions \ref{defnpij} and \ref{defnpin}. \endproof 
 
 The relation between $\pii_n$ and $\pij_n$ is given by 
  \begin{equation}\label{defnpinbis}
 	\pij_n:=\ell \circ \pii_n
 \end{equation}
 \ie composition with the functor  \[\ell:\set2\longrightarrow \Se, \ \ \ell(X):=\varinjlim_{\cC^{o}} X(c)\] which assigns to a $2$-set its set of components. Here, $\cC$ is any of the small categories  defining $\set2$ by duality as $\hat \cC$ as in Proposition \ref{propquivers} and $\cC^{o}$ is its opposite. Note that the functor $\ell$  does not correspond to a point of the topos $\set2$.

\subsection{$\Gamma$-sets as endofuntors } \label{sectendofunct}

In this section we recall the construction of \cite{DGM} of the endofunctor in the category $\sses$  associated to a $\Gamma$-space, in the case of discrete $\Gamma$-spaces, \ie $\sss$-modules. By construction an $\sss$-module is a covariant functor $M:\gop\longrightarrow \Se_*$  and, as in Section 2.1.2.1 of \opcit, we view pointed sets as discrete pointed simplicial sets, \ie as constant functors   $\dop\longrightarrow\Ses$. 
\begin{lem}\label{EilMac}
Let $M:\gop\longrightarrow \Ses$ be an $\sss$-module. Then the associated endofunctor of the category $\sses$ of pointed simplicial sets is obtained by composition with $M$ viewed as 	an endofunctor of the category $\Ses$ of pointed sets. 
\end{lem}
\proof One first extends the functor $M:\gop\longrightarrow \Ses$ to an endofunctor  $\Ses\longrightarrow \Ses$ in pointed sets. This is done by taking a colimit on the finite subsets as explained in \S 2.2.1.1 of \cite{DGM}. Then, one applies the technique described in \opcit that uses, for a simplicial set $X=\{[q]\mapsto X_q\}$, the diagonal
$$
M(X):=\{[q]\mapsto M(X_q)_q\}.
$$
 Since by construction $M(X_q)$ is a discrete simplicial set,  it is the same in all degrees so that the index $q$  in $M(X_q)_q$ disappears, thus we  simply write $M(X_q)$. Hence starting with the pointed simplicial set $X:\dop\longrightarrow \Ses$, we obtain a new pointed simplicial set by composition \ie 
\begin{equation}\label{commp}
 X\mapsto M(X)= M\circ X: \dop \longrightarrow \Ses.
\end{equation} 
In summary the result follows from \S 2.2.1.1 of \opcit. \endproof 
The basic example of an $\sss$-module is given in 2.1.2.1 of \cite{DGM} where one associates to an abelian monoid $A$ with a zero element, the functor $M=HA$  
\begin{equation}\label{hadefn}
HA(k_+)=A^k, \  \  Hf:HA(k_+)\to HA(n_+), \ Hf(m)(j):=\sum_{f(\ell)=j} m_\ell
\end{equation}
where $m=(m_1,\ldots ,m_k)\in HA(k_+)$. The zero element of $A$ gives meaning to the empty sum. In the special case when the monoid $A$ is an abelian group, the composition \eqref{commp}, \ie the functor $HA\circ X: \dop \longrightarrow \Ses$ factors through simplicial abelian groups (the functor $HA$ is the composite of a more precise functor $A\mapsto {\rm Ab}A$ to abelian groups  with the forgetful functor from abelian groups to pointed sets, where the base point is the $0$) and always fulfills the Kan extension property. The geometric realization $\vert HA\circ X\vert$ only uses the underlying simplicial pointed set but  the finer structure as a simplicial abelian group, and the Dold-Kan correspondence in the form of Corollary 2.5 of \cite{GJ}, Chapter III, show  that the homotopy groups of the geometric realization $\vert HA\circ X\vert$ are given by the (reduced) homology\footnote{For a pointed simplicial set $(X,*)$ we use the notation  $H_n(X,A)$ for the reduced homology $H_n((X,*),A)$.} of the associated complex of abelian groups, \ie $\pi_n(\vert HA\circ X\vert)=H_n(X,A)$.
This suffices to conclude for instance that  $\vert HA\circ S^n\vert$ is an Eilenberg-MacLane space $K(A,n)$.

\subsection{The homology with coefficients in an $\sss$-module } \label{secthomology}

 In our arithmetic context we are interested in $\sss$-modules $M$  which are no longer of the form $HA$ where $A$ is an abelian group. In a first class of examples $M$ is still of the form $HA$, where $A$ is a monoid. A second class of examples are those  constructed in \cite{CCprel} to specify the geometric structure of $\overline{\Spec\Z}$ at the archimedean place. In all these cases it is no longer true that the composite $M\circ X$ is fibrant, even when the simplicial set $X$ itself is fibrant.  We shall use the  equality $\pi_n(\vert HA\circ X\vert)=H_n(X,A)$ holding for abelian groups as the motivation to extend the definition of the homology of a pointed simplicial set with coefficients in an arbitrary $\sss$-module as follows
\begin{defn}\label{fact1} Let $M$ be an $\sss$-module, and $X$ a pointed simplicial set. For any integer $n\geq 0$ one defines the homology $H_n(X,M)$ as the  $\sss$-module 
\begin{equation}\label{pins1bis}
H_n(X,M)(k_+):=\pij_n ( M\circ(X\wedge k_+) ).
\end{equation} 
 Here, $k_+$ is viewed as a discrete simplicial pointed set \ie constant in all degrees.
\end{defn}
As in Lemma \ref{EilMac}, $M$ is viewed as 	an endofunctor of the category $\Ses$ and $\pij_n$ is defined in Definition \ref{defnpij}. There is in fact a refined version of homology $H_n^{(2)}(X,M)$ using  $\pii_n$ instead of $\pij_n$ but we shall not need it in the present paper.

The following result establishes several basic properties of the new homology.
\begin{prop}\label{propprel}
$(i)$~For any $n\geq 0$, $H_n(X,M)$ is a covariant bifunctor $$H_n:\sses\times \smod\longrightarrow \smod. $$ 
$(ii)$~Let $M_1,M_2$ be $\sss$-modules. One has a natural transformation
$$
H_n(M_1\circ X,M_2)\to H_n(X,M_2\circ M_1) 
$$
which is an isomorphism when evaluated on $1_+$.\newline
$(iii)$~For any pointed simplicial set $X$ one has $H_n(X,\sss)=\pij_n(X)\wedge \sss$.\newline
$(iv)$~For $n\neq m$ : $H_m(S^n,\sss)=\{*\}$ while for $n=m$ one has $H_m(S^n,\sss)=\sss$.\newline
$(v)$~For $n\neq m$ : $H_m(S^n,H\B)=\{*\}$ while for $n=m$ one has $H_m(S^n,H\B)=H\B$.	
\end{prop}
\proof $(i)$~By construction, $H_n(X,M)$ is a covariant functor of $X$ for fixed $M$, and of $M$ for fixed $X$. To prove that it is a bifunctor it suffices, using the bifunctor lemma (see \cite{ML} Proposition 1 Chapter II, \S 3), to show that it satisfies the interchange law which states that given morphisms $f \in \Hom_\sses(X,Y)$ and $h\in \Hom_\sss(M,N)$ one has the equality 
\begin{equation}\label{biflem}
H_n(f,N)\circ H_n(X,h)=H_n(Y,h)\circ H_n(f,M)\in \Hom_\sss(H_n(X,M),H_n(Y,N)).
\end{equation}
Both sides of  this formula are $\sss$-modules \ie functors $\gop\longrightarrow \Se_*$ thus  it is enough to check the equality pointwise \ie by evaluating both sides on $k_+$  for fixed $k$. Since $\pij_n:\sses\longrightarrow \Se_*$ is a functor the equality follows provided one shows that the same equality holds if one replaces $H_n(X,M)$ by  $F(X,M):=M\circ X$ which is a separately covariant functor to $\sses$ with arguments in $\sses$ and $\smod$. Again it is enough to check this equality pointwise \ie replacing $\sses$ by $\Se_*$ and $F(X,M):=M\circ X$ by $G(X,M):=M(X)$ which is a separately covariant functor to $\Se_*$ with arguments in $\Se_*$ and $\smod$. Since $M$ and $N$ are endofunctors of
$\Se_*$ and the morphism  $h\in \Hom_\sss(M,N)$ is a natural transformation from $M$ to $N$ one has, for any $f \in \Hom_{\Se_*}(X,Y)$, the equality 
$$
N(f)\circ h_X=h_Y\circ M(f)\in \Hom_{\Se_*}(M(X),N(Y))
$$
which gives \eqref{biflem}.\newline
$(ii)$~As in \cite{DGM}, (2.2.1.2 equation (2.2)), one has natural maps $M(X)\wedge Y\mapsto M(X\wedge Y)$. We apply this with $Y=k_+$ and thus obtain natural maps 
$$
\eta_k:M_1(X)\wedge k_+\to M_1(X\wedge k_+) \qqq k.
$$
This yields a natural morphism 
$$
M_2(\eta_k):M_2(M_1(X)\wedge k_+)\to M_2(M_1(X\wedge k_+))
$$
and by composition with $\pij_n$ one gets the natural transformation 
$$
\pij_n(M_2(\eta_k)):H_n(M_1\circ X,M_2)(k_+)\to H_n(X,M_2\circ M_1)(k_+)
$$
which is, by construction, an isomorphism for $k=1$.\newline
$(iii)$~Since the endofunctor of $\Se_*$ associated to $\sss$ is the identity, the result follows from Proposition \ref{propprepa} $(i)$.

$(iv)$~Using $(iii)$ it is enough to determine $\pij_n(S^m)$. The pointed simplicial set $S^n$ is obtained by collapsing $\partial \Delta[n]$ to a base point. This means that one considers the sub-functor   $[q]\mapsto \partial \Delta([q])\subset \Hom_\Delta([q],[n])$ given by the maps $[q]\to [n]$ which are not surjective and one identifies all the elements of $\partial \Delta([q])$ with the base point. For $h\in \Hom_\Delta([q'],[q])$ one has $\partial \Delta([q])\circ h\subset \partial \Delta([q'])$ so that the collapsing gives a pointed simplicial set. An element of $\Hom_\sses(S^n,X)$ is an element 
of $\Hom_\sses(\Delta[n],X)$ which maps $\partial \Delta[n]$ to the base point. This means, by Yoneda's lemma, an element  $x\in X_n$ such that $\partial_j(x)=*$ for all $j$ (since any map $[q]\to [n]$ which is not surjective factors through a $d_j:[n-1]\to [n]$). For $X=S^m$, such an $x\in X_n$  is, if it is not the base point, an element $\phi \in \Hom_\Delta([n],[m])$ which is surjective and such that $\phi\circ d_j$ fails to be surjective for any $j$. This latter condition implies that $\phi$ is also injective and one concludes that $n=m$ and $\phi$ is the identity map. This gives $\pij_n (S^m)=\{*\}$ for $n\neq m$.  To prove that $\pij_n (S^n)=\{*,\id\}$  one just needs to show that the element $\id$ does not get identified with the base point under the equivalence relation generated by the relation  \eqref{Maydefnbis}, \ie
$$
  xRy\iff  \exists z \mid \partial_jz=*\, \forall j<n, \ \partial_nz=x,\, \partial_{n+1}z=y\}.
 $$ 
 Any $z\neq *$ in $S^n_{n+1}$ is given by a surjective map $s_i \in \Hom_\Delta([n+1],[n])$ such that $s_i (i)=s_i(i+1)$ and the condition $\partial_jz=*\, \forall j<n$ shows that the index $i$ is equal to $i=n$. It follows that $\partial_nz=\partial_{n+1}z$ and that the relation $R$ is the diagonal. \newline
 $(v)$~The endofunctor $H\B$ associates to a pointed set $E$ the (pointed) set of all finite subsets of $E$ which contain the base point $*$, and to a map $f:E\to F$ the direct image map $Z\mapsto f(Z)$. Note the equivalence, 
 \begin{equation}\label{charoneequi}
 	f(Z)=\{*\}\iff f(x)=*\qqq x\in Z
 \end{equation}
 It follows that there are only two elements $u\in (H\B\circ S^n)_n=H\B(S^n_n)=H\B(\{*,\id\})$, namely the base point $*$ and the subset $u=\{*,\id\}$.  
Let us show that these two elements are not equivalent under the equivalence relation generated by the relation  \eqref{Maydefnbis} for the simplicial pointed set  $H\B\circ S^n$. An element of $(H\B\circ S^n)_{n+1}=H\B(S^n_{n+1})$ is a subset $z=\{*,s_{i_1}, \ldots , s_{i_k}\}$ of the set $S^n_{n+1}$ described in the proof of $(iv)$. The condition 
$\partial_jz=*\, \forall j<n$ shows that all indices $i_j$ are equal to $n$ so that either $z=*$ or $z=\{*,s_n\}$. This shows, as in the proof of $(iv)$, that the relation $R$ is diagonal and  $\pij_n(H\B\circ S^n)=H\B(1_+)=\B$. Let then $k>0$ be an integer and $E$ a pointed set. One has a natural isomorphism $H\B\circ (E \wedge k_+)\simeq (H\B\circ E)^k$. It follows that for any pointed simplicial set $X$ one has a natural isomorphism $H\B\circ (X\wedge k_+)\to (H\B\circ X)^k$. Then by Proposition \ref{propprepa} $(ii)$ one gets 
$$
 \pij_n(H\B\circ (S^n \wedge k_+))=(\pij_n(H\B\circ S^n))^k=H\B(k_+). 
$$
By construction the natural identifications are compatible with the structures of $\Gamma$-sets. This shows that $H_n(S^n,H\B)=H\B$. The proof of $(iv)$ together with \eqref{charoneequi} show that $H_n(S^n,H\B)=\{*\}$ for $m\neq n$. \endproof 

 Definition \ref{fact1} provides a  meaning to the following equality \eqref{equbasic} whose two sides  are $\sss$-modules.
\begin{thm}\label{propcompare}
	Let $A$ be an abelian group, and $X$ a pointed simplicial set. For any integer $n\geq 0$ one has the equality of $\sss$-modules
	\begin{equation}\label{equbasic}
  H_n(X,HA)=H(H_n(X,A))
\end{equation} 
where $H_n(X,A)$ is the (reduced) abelian group homology  of $X$ with coefficients in $A$.
\end{thm}
\proof For any simplicial set $Y$ the composite $HA\circ Y$ is a simplicial abelian group and hence has the Kan extension property (see \cite{Moore}, Theorem 2.2). It follows that the combinatorial homotopy $\pij_n(HA\circ Y)$ coincides with the usual homotopy $\pi_n(\vert HA\circ Y\vert)$ of the geometric realization
\begin{equation}\label{usualhomot}
	\pij_n(HA\circ Y)=\pi_n(\vert HA\circ Y\vert).
\end{equation}
 Moreover   the group law of these homotopy groups coincides with the abelian group law inherited from the simplicial abelian group structure (see \opcit Proposition 2.4). The Dold-Kan correspondence   (see \cite{GJ}, Chapter III Corollary 2.5) gives a canonical bijection 
 \begin{equation}\label{usualhomot1}
\delta_Y: H_n(Y,A)\to \pij_n(HA\circ Y). 
\end{equation}
Furthermore this bijection is a natural transformation of covariant functors from pointed simplicial sets to pointed sets. More precisely given a morphism $\psi:Y\to Y'$ of pointed simplicial sets, 
one obtains the equality 
\begin{equation}\label{usualhomot2}
\pij_n(HA(\psi))\circ \delta_Y=\delta_{Y'}\circ H_n(\psi,A).
\end{equation}
Indeed, it is enough to check this equality on cycles $c\in Z_n(Y,A)$ which are Moore normalized, \ie $\partial_j c=0$ $\forall j$. The element $\delta_Y(c)$ is then given by the combinatorial class directly associated to $c$ viewed as an element of $\Hom_\sses(S^n,HA\circ Y)=\Hom_\sses((\Delta^n,\partial \Delta^n),HA\circ Y)$ which is a subset of $(HA\circ Y)_n=HA(Y_n)$. Let $c=\sum a_j y_j$ where $a_j \in A$ and $y_j\in Y_n$. Then $\pij_n(HA(\psi))\circ \delta_Y(c)$ is represented by the combinatorial class obtained from $\delta_Y(c)$ by applying the functor $HA(\psi)$. This gives $HA(\psi)(\sum a_j y_j)=\sum a_j \psi(y_j)$, as one sees using the definition \eqref{hadefn} of the functor $HA$. But one has similarly $H_n(\psi,A)(c)=\sum a_j \psi(y_j)$, thus one gets the required equality \eqref{usualhomot2}.

From \eqref{usualhomot1} one gets the bijection  $\delta_X: H_n(X,A)\to \pij_n(HA\circ X)$. Let then $k>0$ be an integer and $E$ a pointed set. One has a natural isomorphism $HA\circ (E \wedge k_+)\simeq (HA\circ E)^k$ since both sides consist of maps $(x,j)\mapsto \phi(x,j)\in A$, $x\in E$, $j\in \{1, \ldots, k\}$ with finite support and such that $\phi(*,j)=0$ for all $j$. Thus one obtains a natural isomorphism of simplicial sets
$$
HA\circ (X\wedge k_+)=(HA\circ X)^k.
$$
The same equality holds for the geometric realizations, and using \eqref{usualhomot} one derives
$$
\pij_n(HA\circ (X\wedge k_+))=\pi_n(\vert HA\circ (X\wedge k_+)
\vert)= \pi_n((\vert HA\circ X\vert)^k)= H_n(X,A)^k.
$$
 At the set-theoretic level this coincides with $H(H_n(X,A))(k_+)$. In fact one can obtain the same result more directly as a consequence of \eqref{usualhomot1} and of the equality of (reduced) homology groups $H_n(X_1\vee X_2,A)=H_n(X_1,A)\oplus H_n(X_2,A)$.\newline
  It remains to show that given a morphism $\phi:k_+\to m_+$ in $\gop$ the associated map  
  $$\pij_n(HA\circ (X\wedge k_+))\to \pij_n(HA\circ (X\wedge m_+)),
  $$ 
  is the same as the map $HK(\phi): H_n(X,A)^k\to H_n(X,A)^m$ associated to the group law of $K=H_n(X,A)$ and the functor $HK$. Using \eqref{usualhomot2} it is enough to show that $HK(\phi)$  equals the homology map 
 $$
 H_n(\id_X \wedge \phi,A): H_n(X\wedge k_+,A)\to H_n(X\wedge m_+,A).
 $$
 With $e_j$, $j\in \{1,\ldots, k\}$, the canonical basis of $H_0(k_+)$, and $e_0:=0$, the above map is given by 
 $$
 H_n(\id_X \wedge \phi,A)(\sum c_j\otimes e_j)=\sum c_j\otimes e_{\phi(j)}
 $$
 and using the definition \eqref{hadefn} of the functor $HK$ one gets the required equality.\endproof 
 \begin{rem} \label{homotopgroup} In homotopy theory the homotopy groups $\pi_n$  are abelian groups for $n>1$. The group operation arises, at the combinatorial level, from the Kan extension property of fibrant simplicial sets together with combinatorial constructions involving simplices. Definition \ref{fact1} does not involve any of these constructions and yet Theorem \ref{propcompare} shows that one recovers the same group law on the homotopy groups $\pi_n$  from the  $\Gamma$-set ($\sss$-module) obtained using the functorial nature of the map $k_+\mapsto X\wedge k_+$. The reason behind this equality of structures is the interchange law which is fulfilled by the group law of the homotopy group $\pi_n$ and the group law induced by the abelian coefficients. In that sense, Definition \ref{fact1} takes into account the $\sss$-module structure of the coefficients to obtain a replacement of the group structure of homotopy groups. We shall see in the next sections a striking example where this additional structure is put to work. 	
 \end{rem}

\section{The archimedean place and the Gromov norm}\label{sectgromov}
In this section we show that the singular homology $H_*(X, \R)$ of a topological space inherits a natural semi-norm from the filtration of the $\sss$-module $H\R$ by the sub-$\sss$-modules $\Vert H\R \Vert_\lambda$ associated to the archimedean place of $\overline{\Spec\Z}$ as constructed in \cite{CCprel}. Moreover we prove that this semi-norm is equivalent to the Gromov semi-norm on singular homology.

\subsection{$\overline{\Spec\Z}$ at the archimedean place}\label{sectspecZ}
In \cite{CCprel} we showed how to endow the Arakelov compactification $\overline{\Spec\Z}$ with a structure sheaf of $\sss$-algebras, which coincides with the standard structure sheaf of $\Spec\Z$ on the dense open set $\Spec\Z\subset\overline{\Spec\Z}$ using the fully faithful functor $H$ from rings to $\sss$-algebras. The new feature is the structure of this sheaf at the archimedean place which is obtained using the following  proposition\footnote{With the nuance that  in \eqref{subvert1} we use the strict inequality.} of \cite{CCprel}  
\begin{prop}\label{sssalg2} $(i)$~Let $R$ be a semiring, and $\vvert\ \ \vvert$ a sub-multiplicative seminorm on $R$. Then $HR$ is naturally endowed with a structure of $\sss$-subalgebra $\vvert HR\vvert_1\subset HR$ defined as follows
\begin{equation}\label{subvert}
\vvert HR\vvert_1: \Gamma^{{\rm op}}\longrightarrow \Se_*\qquad \vvert HR\vvert_1(X):=\{\phi\in HR(X)\mid \sum_{X\setminus \{*\}} \vvert\phi(x)\vvert\leq 1\}.
\end{equation}
$(ii)$~Let $E$ be an $R$-semimodule and $\vvert\ \ \vvert^E$ a seminorm on $E$ such that $\vvert a\xi\vvert\leq \vvert a\vvert \vvert \xi\vvert$, $\forall a\in R$, $\forall \xi \in E$, then for any $\lambda\in \R_+$ the following defines a module $\vvert HE\vvert^E_\lambda$ over $\vvert HR\vvert_1$
\begin{equation}\label{subvert1}
\vvert HE\vvert^E_\lambda: \Gamma^{{\rm op}}\longrightarrow \Se_*\qquad \vvert HE\vvert^E_\lambda(X):=\{\phi\in HE(X)\mid \sum_{X\setminus \{*\}} \vvert\phi(x)\vvert^E< \lambda\}.
\end{equation}
\end{prop}
The first statement of  Proposition \ref{sssalg2}  is applied for the ring $R=\Q$ of rational numbers and its archimedean absolute value to construct the stalk at $\infty$ of the structure sheaf. One  obtains in this way  a sheaf $\hq$ of $\sss$-algebras over $\overline{\Spec\Z}$.
 The second statement of  Proposition \ref{sssalg2}  is then applied to the one-dimensional real vector space  $\R$ to obtain, given an Arakelov divisor $D=D_\fii+D_\infty$, the  sheaf $\cO(D)$ of $\cO$-modules over $\overline{\Spec\Z}$
\begin{equation}\label{sheafsmod}
\cO(D)(\Omega):=\vvert H\cO(D_\fii)(\Omega\setminus\{\infty\})\vvert_{e^a}, \  \ D_\infty=a\{\infty\}.
 \end{equation}
Thus the $\sss$-modules at work at the archimedean place depend on a positive real parameter $\lambda>0$ and are implemented  by the 
functor $\Vert H\R \Vert_\lambda:\gop \longrightarrow \Se_*$ which associates to a pointed set $X$ the pointed set 
\begin{equation}\label{hrlam}
\Vert H\R \Vert_\lambda(X)=\{x:X\to \R\mid \#\{j,x(j)\neq 0\}<\infty\ \& \sum \vert x(j)\vert < \lambda\}.  
 \end{equation}
In fact \eqref{hrlam} describes also the extension of $\Vert H\R \Vert_\lambda$ as an endofunctor of $\Se_*$. There is an obvious analogue of  \eqref{hrlam} when $\R$ is replaced by $\Q$ and this analogue is what is needed in \eqref{sheafsmod}; on the other hand it is more natural to work  with the local field $\R$ associated to the archimedean place of $\Q$.

\subsection{Simplicial volume}\label{simplicialvol}
We recall briefly the notion of simplicial volume introduced by M. Gromov \cite{gromov82}. Let first $X$ be a topological space and $C_*(X,\R)$ the associated singular chain complex with real coefficients. One defines the $\ell^1$-norm on singular chains as follows
\begin{equation}\label{ell1}
\Vert c\Vert_1:=\sum \vert a_j\vert \qqq c= \sum a_j \sigma_j, \ \sigma_j\in {\rm Top}(\Delta^*,X)
\end{equation}
The induced semi-norm on the singular homology $H_*(X, \R)$ is the quotient semi-norm
\begin{equation}\label{ell1bis}
\Vert \alpha\Vert_1:=\inf_{c\in \alpha} \Vert c\Vert_1.
\end{equation}
The Gromov norm $\vert M\vert$ of an oriented closed connected manifold of dimension $n$ is then defined as the semi-norm of its fundamental class $\vert M\vert:=\Vert [M]\Vert_1$.
A fundamental result of the theory (\!\cite{gromov82}, \cite{thurston} Thm 6.2) is the proportionality principle:
 \begin{thm}\label{gromov} (M. Gromov) Let $\Sigma$ be any compact oriented hyperbolic manifold of dimension $n>1$, then one has 
  $$\vert \Sigma\vert= \frac{v(\Sigma)}{v_n}$$
 	where $v(\Sigma)$ is the volume of $\Sigma$ and $v_n$ is the maximal volume of straight simplices in hyperbolic space. 
 \end{thm}
 We refer to \cite{thurston} chapter 6 for the description of the straightening of singular simplices and singular chains. The constant $v_2$ is equal to $\pi$ and one thus has 
 \begin{cor}\label{gromov1}
 	Let $\Sigma$ be a Riemann surface of genus $g>1$, then $\vert \Sigma\vert=4(g-1)$.
 \end{cor}
The fact that the norm does not vanish is dual to the boundedness of cohomology and this holds  in the hyperbolic case, thus for $k>1$ the semi-norm  \eqref{ell1bis} is in fact a norm on the homology $H_k(M,\R)$ when $M$ is an hyperbolic manifold (see \opcit). 
 
 \subsection{Moore normalization}\label{sectmoore}
 
 Let $A$ be a simplicial abelian group. The standard complex (still denoted $A$ for simplicity) of abelian groups associated to $A$ is defined using the boundary map
 \begin{equation}\label{standard}
\partial:=\sum_{0}^{ n}\,(-1)^j d_j: A_n \to A_{n-1}.
\end{equation}

  The associated normalized complex $NA$ is defined as follows
\begin{equation}\label{normal}
NA_n:=\cap_{0}^{ n-1} \,\ker\, d_j\subset A_n,\ \  d:=d_n:NA_n\to NA_{n-1}
\end{equation}
(the simplicial identity $d_{n} d_{n-1}=d_{n-1}d_{n-1}$ shows that it defines a complex).
 For each $n$ one lets $D_n\subset A_n$ be the subgroup generated by the ranges of the degeneracies. The boundary map $\partial$ of \eqref{standard} fulfills $ \partial(D_n)\subset  D_{n-1}$ and induces a map
 $$
 \partial: A_n/D_n\to A_{n-1}/D_{n-1}.
 $$
 The corresponding quotient complex $A/D$ is the complex modulo degeneracies. By construction one has two morphisms of complexes $i:NA\to A$ and $p:A\to A/D$. Moreover (see \cite{GJ}, Theorem 2.1), the morphism $p\circ i: NA\to A/D$ is an isomorphism of chain complexes. It follows that  the composite morphism 
 $\nu:=(p\circ i)^{-1}\circ p:A\to NA$ is a projection. As in  the proof of Theorem 2.4 of \cite{GJ}, one constructs explicitly a chain map $f^{(n)}:A_n\to NA_n$ as the composition
 \begin{equation}\label{normal1}
 f^{(n)}:=f^{(n)}_{n-2}\circ \ldots \circ f^{(n)}_j\circ f^{(n)}_0
 \end{equation}
 where $f^{(n)}_j:A_n\to A_n$ is defined as $f^{(n)}_j=\id -s_{j+1}d_{j+1}$. Moreover one also constructs explicitly a chain homotopy $T_k:A_k\to A_{k+1}$ such that 
 \begin{equation}\label{normal2}
 \id-f^{(n)}=T\circ \partial +\partial\circ T.
 \end{equation}
 Since each $f^{(n)}_j$ acts as the identity in the quotient $A_n/D_n$ the same holds for $f^{(n)}$ and one obtains the equality $\nu_n=f^{(n)}$. 
\begin{lem}\label{seminorm} Let $X$ be pointed simplicial set. Let the simplicial vector space $A=H\R\circ X$ be endowed with the norm 
\begin{equation}\label{norm1}
\Vert \phi\Vert := \sum \vert \phi(x)\vert \qqq \phi \in A_n=H\R(X_n).
\end{equation}
$(i)$~The linear map $\nu_n=f^{(n)}:A_n\to NA_n$ is of norm $\leq 2^{n-1}$.\newline
$(ii)$~Let $c\in Z_n(A)$ be a cycle. Then $\nu_n(c)$ is a homologous normalized cycle and $\Vert \nu_n(c)\Vert\leq 2^{n-1}\Vert c\Vert $.
\end{lem}
\proof $(i)$~The statement follows from \eqref{normal1} and the inequalities
$$
\Vert  f^{(n)}_j\Vert \leq \Vert  \id\Vert+\Vert  s_{j+1}d_{j+1}\Vert\leq 2, \ \ \Vert  f\circ g\Vert\leq \Vert f\Vert\Vert  g\Vert.
$$ 
$(ii)$~This follows from $(i)$ and \eqref{normal2}. \endproof 
On the real vector space $H_n(X,\R)$ we consider the following semi-norm which is induced by the $\ell^1$-norm \eqref{norm1} on the normalized complex $NA_n$:
\begin{equation}\label{normN}
\Vert c\Vert^{\rm nor}:=\inf \{\Vert \phi\Vert \mid \phi \in NA_n, \ \phi \sim c\}
\end{equation}
where $\phi \sim c$ means that $\phi \in NA_n$ is homologous to the cycle $c\in Z_n(A)$. By applying Lemma \ref{seminorm} one obtains the basic inequalities 
\begin{equation}\label{normN1}
\Vert c\Vert_1\leq \Vert c\Vert^{\rm nor}\leq 2^{n-1} \Vert c\Vert_1\qqq c\in H_n(X,\R).
\end{equation}

\subsection{Equivalence with the Gromov norm}\label{sectequivalence}
 The filtration of the $\sss$-module $H\R$ by the sub-$\sss$-modules 
 $\Vert H\R \Vert_\lambda\subset  H\R $, $\lambda\in \R_+$, of \eqref{hrlam}  provides, for any pointed simplicial set $X$ and integer $n\geq 0$ natural morphisms of $\sss$-modules  
\begin{equation}\label{filtr1}
\rho_{n,\lambda}:H_n(X,\Vert H\R \Vert_\lambda)\to H_n(X, H\R), 
\end{equation}
and a filtration by the ranges of the $\rho_{n,\lambda}$. Theorem \ref{propcompare} gives  a natural isomorphism $ H_n(X,H\R)=H(H_n(X,\R))$. 
\begin{thm}\label{propcompare1}
	Let  $X$ be a pointed simplicial set. For any integer $n\geq 0$ and $\lambda\in \R_+$, the range of the natural morphism of $\sss$-module is
		\begin{equation}\label{equbasic2}
 \rho_{n,\lambda}\left(H_n(X,\Vert H\R \Vert_\lambda)\right) =\Vert H(H_n(X,\R))\Vert_\lambda^{\rm nor}
\end{equation} 
where $\Vert c\Vert^{\rm nor}$  is the semi-norm defined in \eqref{normN}.
\end{thm}
\proof 
It follows from \eqref{hrlam} that the endofunctor $\Vert H\R \Vert_\lambda$  assigns to a pointed set $(X,*)$ the set of maps with finite support
$$
\Vert H\R \Vert_\lambda(X):=\{\phi:X\to \R\mid \phi(*)=0, \ \#\{x\mid \phi(x)\neq 0\}<\infty, \  \sum \vert \phi(x)\vert < \lambda\}.
$$
By construction the range $\rho_{n,\lambda}\left(H_n(X,\Vert H\R \Vert_\lambda)\right)$ is a sub-functor of $H(H_n(X,\R))$ thus to show \eqref{equbasic2} it is enough to prove that for any integer $k>0$ one has
\begin{equation}\label{equbasic3}
 \rho_{n,\lambda}\left(H_n(X,\Vert H\R \Vert_\lambda)\right)(k_+) =\Vert H(H_n(X,\R))\Vert_\lambda^{\rm nor}(k_+).
 \end{equation} 
The right hand side of \eqref{equbasic3} is given by $k$-tuples $(\gamma_j)_{1\leq j\leq k}$, $\gamma_j\in H_n(X,\R)$ such that $\sum \Vert \gamma_j\Vert^{\rm nor}< \lambda$, \ie 
$$
 \exists \phi_j\in H\R\circ X_n\mid d_i\phi_j=0 \, \forall i, \ \phi_j\sim \gamma_j, \ \sum_{X_n\times \{1,\ldots,k\}} \vert \phi_j(x)\vert < \lambda.
$$
For the left hand side of \eqref{equbasic3} one has 
$$
\Vert H\R \Vert_\lambda(X_n\wedge k_+)=\{(\psi_j)_{j\in \{1,\ldots,k\} }\mid \psi_j:X_n\to \R, \, \psi_j(*)=0, \ \#\{x\mid \psi_j(x)\neq 0\}<\infty, \  \sum \vert \psi_j(x)\vert < \lambda\}
$$ 
and moreover the simplicial structure satisfies 
$$
d_i\left((\psi_j)_{j\in \{1,\ldots,k\} }\right)=(d_i\psi_j)_{j\in \{1,\ldots,k\}}.
$$
Thus  the $0$-chains
$
\Hom_\sses(S^n,\Vert H\R \Vert_\lambda(X\wedge k_+) )
$
are exactly the same as the ones involved in the right hand side of \eqref{equbasic3} and one gets  \eqref{equbasic2}.
\endproof

For a topological space $X$ one lets $\sin X$ be the associated simplicial set of singular simplices  
$$
 \sin X=\{[n]\mapsto {\rm Top}(\Delta^n,X)\}
 $$
 where the standard simplex $\Delta^n$ of dimension $n$ is given concretely as 
 $$
 \Delta^n=\{(\lambda_0,\ldots,\lambda_n)\mid \lambda_j\geq 0 , \ \sum \lambda_j=1\}
 $$ 
Then, Definition \ref{fact1} extends to topological spaces and arbitrary $\sss$-modules as
\begin{equation}\label{pins2}
H_n(X,M):=H_n (\sin  X,M ) 
\end{equation}

\begin{cor}\label{conj} Let $X$ be a topological space. The filtration  of the singular homology group $H_n(X,\R)$ by the $\sss$-modules $H_n(X,\Vert H\R \Vert_\lambda)$ defines a semi-norm 
which is equivalent to the Gromov norm. 
 \end{cor}
\proof This follows from Theorem \ref{propcompare1} and the basic inequalities \eqref{normN1}.\endproof 

\vspace{0.1cm}

\section{Equality for Riemann surfaces of genus $g>1$}\label{sectequal}
 We show that for a compact oriented $2$-dimensional manifold $\Sigma$ of genus $g>1$, the normalized norm \eqref{normN} on singular homology agrees with the Gromov norm. Since these two norms are equivalent and the Gromov norm vanishes except on $H_2(\Sigma,\R)$ it is enough to prove the equality for the fundamental class $[\Sigma]\in H_2(\Sigma,\R)$. The difficulty is to construct singular cycles $c$ in the homology class $[\Sigma]$ which not only have $\ell^1$-norm $\Vert c\Vert_1$ close to the expected value $4(g-1)$ but are also {\em normalized}, \ie such that all boundaries vanish $\partial_j(c)=0$. This is achieved in three steps. In section \ref{block} we deal with the relative situation of the building block $K$ and construct a normalized cycle relative to    its boundary $\partial K$. In section \ref{block1} we assemble together $g$ copies of  $K$ and obtain  a surface of genus $g$ and a normalized cycle representing the fundamental class whose $\ell^1$-norm is $4g$. The third step  is standard and uses cyclic covers to improve the estimate to the expected value $4(g-1)$.

\subsection{Moore normalization for the building block}\label{block}
\begin{figure}[!htb]
   \begin{minipage}{0.48\textwidth}
     \centering
     \includegraphics[width=.6\linewidth]{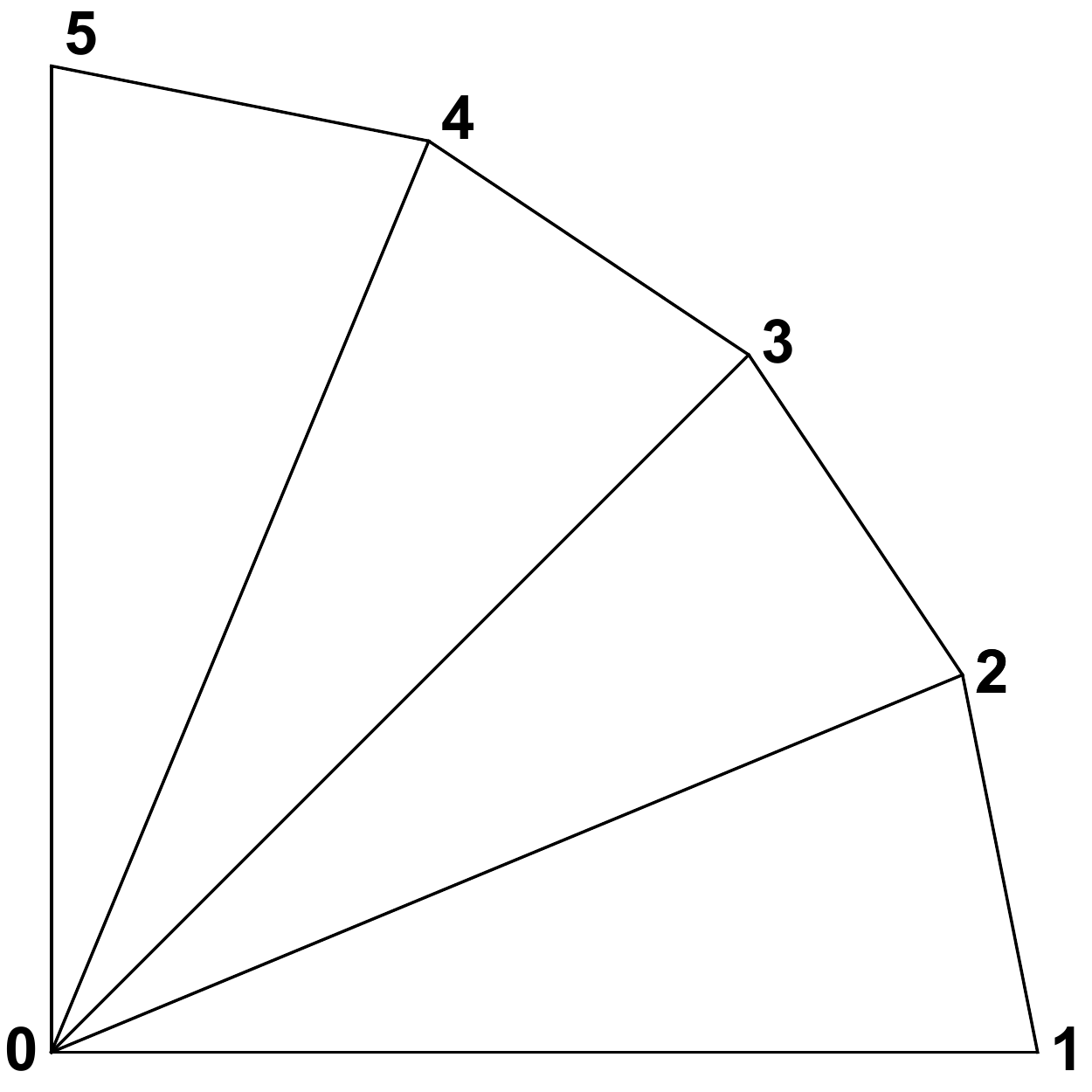}
     \caption{Basic polygon in $E=\R^2$.}\label{gromovnorm1}
   \end{minipage}\hfill
   \begin{minipage}{0.48\textwidth}
     \centering
     \includegraphics[width=.6\linewidth]{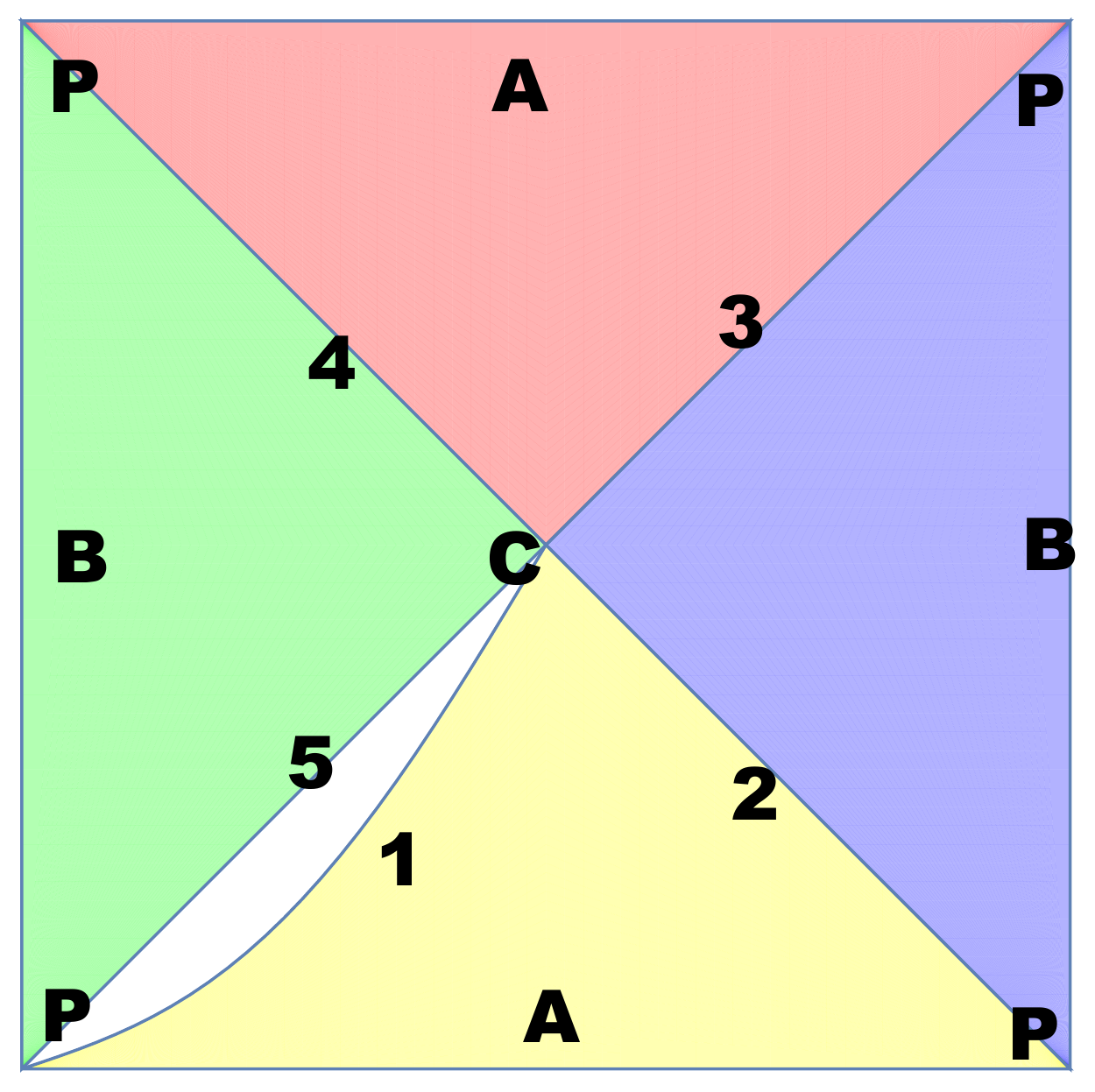}
     \caption{Basic polygon inside the fundamental domain for $\Z^2$ acting on $\R^2$.}\label{flattriang1}
   \end{minipage}
\end{figure}

\begin{figure}[t!]
   \begin{minipage}{0.48\textwidth}
     \centering
     \includegraphics[width=.8\linewidth]{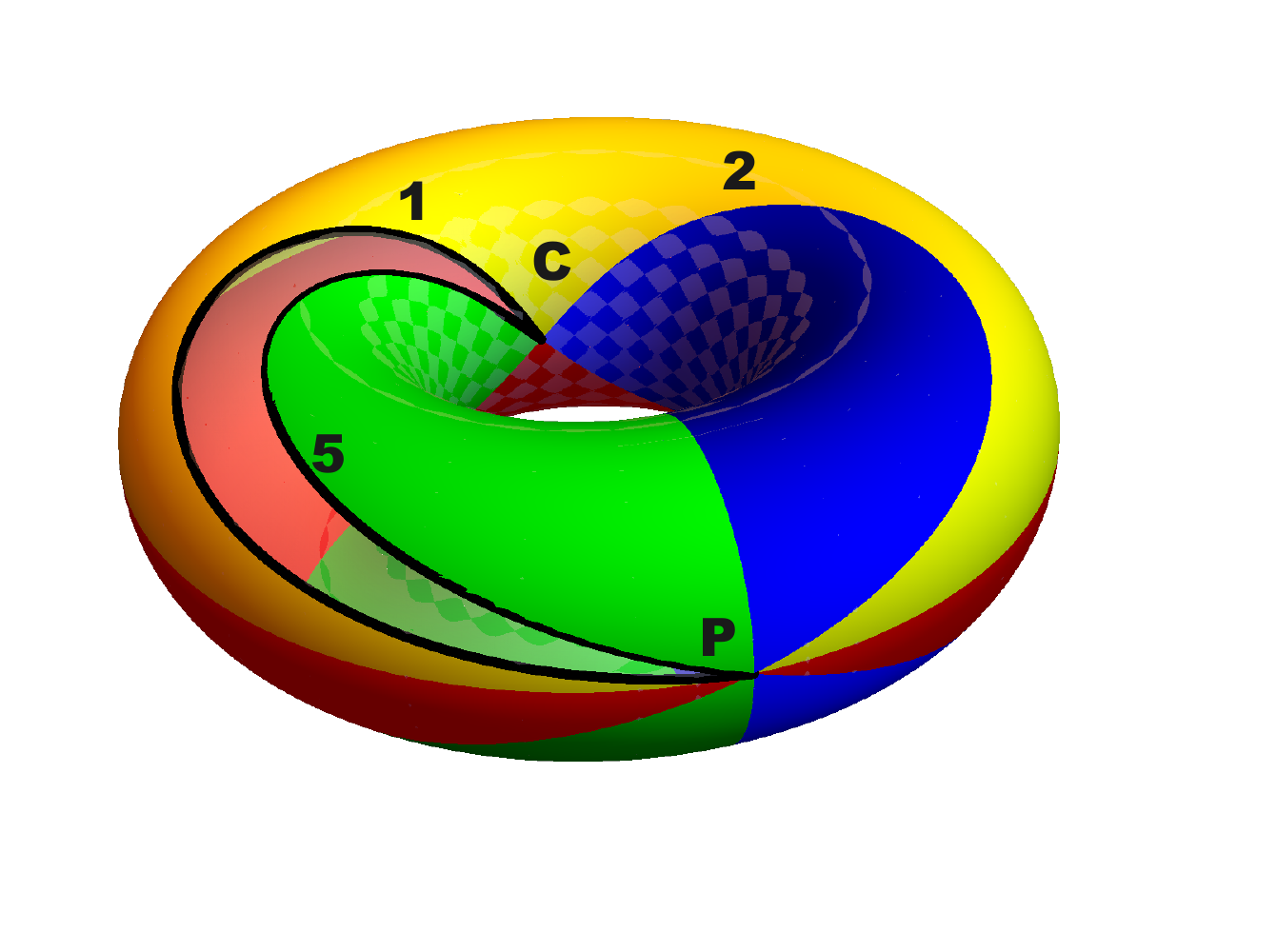}
     \caption{Triangulated building block $K$.}\label{cuttorus1}
   \end{minipage}\hfill
   \begin{minipage}{0.48\textwidth}
     \centering
     \includegraphics[width=.9\linewidth]{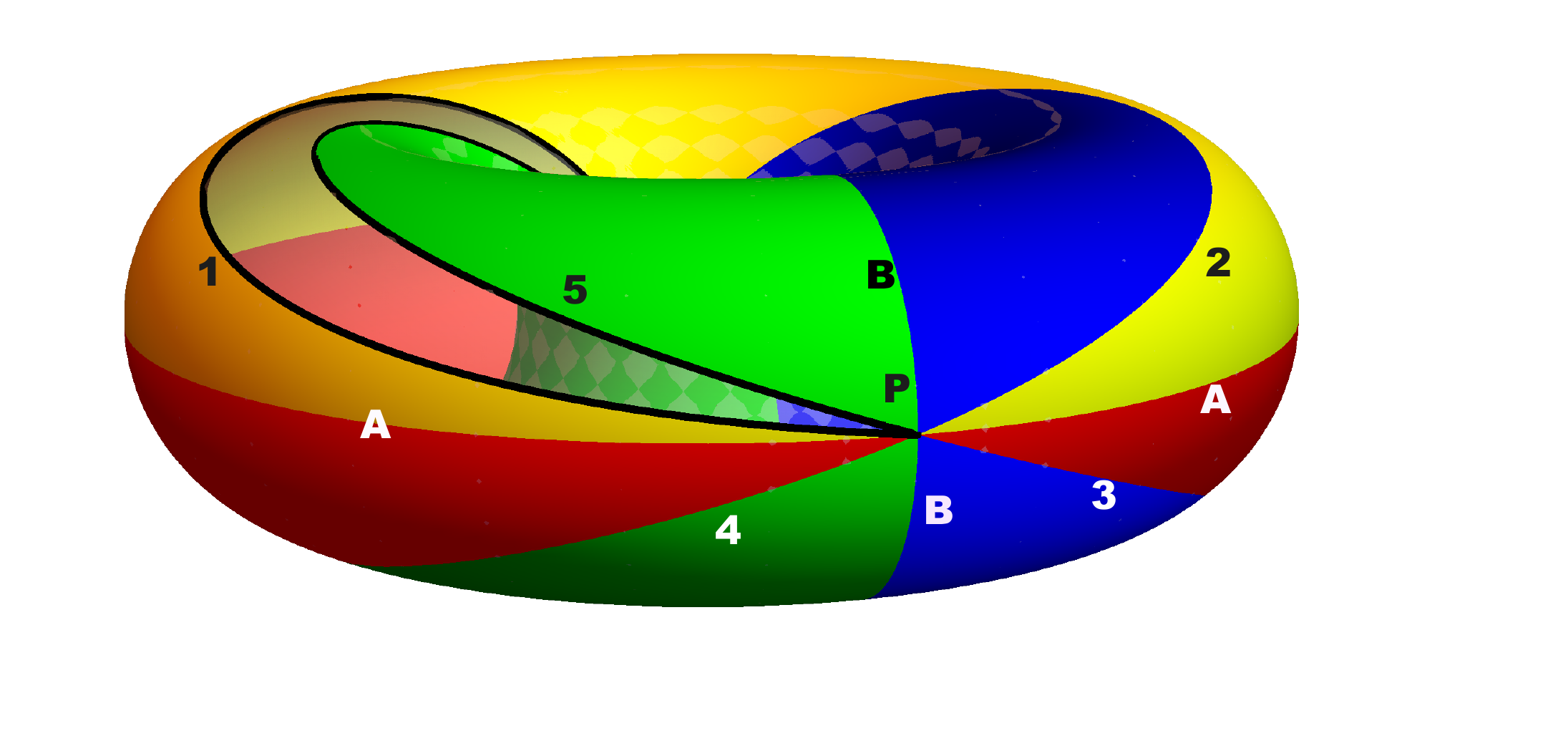}
     \caption{Neighborhood of point $P$.}\label{triangulationtor}
   \end{minipage}
\end{figure}

A compact oriented $2$-dimensional manifold $\Sigma$ of genus $g > 1$ is obtained by gluing together $g$ copies of a building block $K$ which we now describe.  
 This building block is the quotient of the convex polygon 
 ${\rm Conv}(0,1,2,3,4,5)$ of Figure \ref{gromovnorm1} by the equivalence relation $R$ generated by 
 $$
 \Delta(\{1,2\})(x)\sim_R 	\Delta(\{4,3\})(x), \ \ \Delta(\{2,3\})(x)\sim_R\Delta(\{5,4\})(x) \qqq x\in \Delta^1
 $$
 where given $n+1$ points $(P_0,\ldots,P_n)$ in the real affine plane $E=\R^2$, one denotes 
 $$
 \Delta(\{P_0,\ldots,P_n\})\in {\rm Top}(\Delta^n,E), \quad (\lambda_0,\ldots,\lambda_n)\mapsto \sum \lambda_j P_j.
 $$  
  By transitivity one finds that the five vertices $(1,2,3,4,5)$ are equal modulo $R$, since $1 \sim_R 4 \sim_R 3\sim_R 2\sim_R 5$. Thus one has by constrution a continuous map 
 \begin{equation}\label{gamma}
 	\gamma: {\rm Conv}(0,1,2,3,4,5)\to K.
 \end{equation}
  The building block $K$ thus consists of $4$ triangles with the common vertex $0$ and where the external sides are identified following the rules 
\begin{equation}\label{rule}
\Delta(\{1,2\})\sim 	\Delta(\{4,3\}), \ \Delta(\{2,1\})\sim\Delta(\{3,4\}), \ \Delta(\{2,3\})\sim\Delta(\{5,4\}), \ \Delta(\{3,2\})\sim\Delta(\{4,5\}).
\end{equation}
It is 
 shown geometrically  in Figure \ref{flattriang1} as a subset of the $2$-torus (before the identifications of the edges) and in Figures \ref{cuttorus1} and \ref{triangulationtor} after these identifications have been performed.  These Figures keep track of the natural  triangulations.

   By construction one has $\Delta(\{P_0,\ldots,P_n\})\in {\rm Top}(\Delta^n,{\rm Conv}(P_j))$ where ${\rm Conv}(P_j)$ is the convex hull of the points $P_j$. The composition 
 \begin{equation}\label{ruleprime}
\dellt(\{P_0,\ldots,P_n\}):=\gamma\circ \Delta(\{P_0,\ldots,P_n\})
\in {\rm Top}(\Delta^n,K) 
\end{equation}
defines  singular simplices \ie  elements of $\sin\, K$. From \eqref{rule} one obtains the equalities
\begin{equation}\label{rulepri}
\dellt(\{1,2\})=	\dellt(\{4,3\}), \ \dellt(\{2,1\})=\dellt(\{3,4\}), \ \dellt(\{2,3\})=\dellt(\{5,4\}), \ \dellt(\{3,2\})=\dellt(\{4,5\})
\end{equation}

To each pair  $(i,j)$ of $>0$ integer indices  we associate the simplicial chain $c(0,i,j)$
	\begin{equation}\label{chaingood}
c(0,i,j):=	\dellt (\{0,i,j\})+\dellt (\{i,j,0\})+2 \dellt (\{j,0,i\})-\dellt (\{j,i,0\})-\dellt (\{0,j,i\})-2 \dellt (\{i,0,j\})
	\end{equation}

 The boundaries of $c(0,i,j)$ are described as follows

\begin{lem}\label{chaingood1} The following equalities  hold 
\begin{align*}
\partial_0 (c(0,i,j))&=2 \dellt (\{0,i\})-2 \dellt (\{0,j\})-\dellt (\{i,0\})+\dellt (\{j,0\})+\dellt (\{i,j\})-\dellt (\{j,i\})
\\
\partial_1 (c(0,i,j))&=-\dellt (\{0,i\})+\dellt (\{0,j\})+\dellt (\{i,0\})-\dellt (\{j,0\})+2 \dellt (\{j,i\})-2 \dellt (\{i,j\})
\\
\partial_2 (c(0,i,j))&=\dellt (\{0,i\})-\dellt (\{0,j\})+2 \dellt (\{j,0\})-2 \dellt (\{i,0\})+\dellt (\{i,j\})-\dellt (\{j,i\})
\end{align*}
\end{lem}
\proof The result follows by linearity of the $\partial_j$ and the equalities
$$
\partial_0 \dellt (\{a,b,c\})=\dellt (\{b,c\}), \ \partial_1 \dellt (\{a,b,c\})=\dellt (\{a,c\}), \ \partial_2 \dellt (\{a,b,c\})=\dellt (\{a,b\})
$$
\endproof 
We now combine the above simplicial chains and use the rules \eqref{rulepri} to get a chain which is normalized relative to the boundary $\partial K$  of $K$.
 \begin{lem}\label{chaingood2} Let $c_{(1,5)}:=\sum_{1\leq i\leq 4} c(0,i,i+1)$. One has
 \begin{align*}
 \partial_0 c_{(1,5)}&=2 \dellt (\{0,1\})-2 \dellt (\{0,5\})-\dellt (\{1,0\})+\dellt (\{5,0\})
 \\
 \partial_1 c_{(1,5)}&=-\dellt (\{0,1\})+\dellt (\{0,5\})+\dellt (\{1,0\})-\dellt (\{5,0\})
 \\
 \partial_2  c_{(1,5)}&=
\dellt (\{0,1\})-\dellt (\{0,5\})-2 \dellt (\{1,0\})+2 \dellt (\{5,0\})
 \end{align*} 	 \end{lem}
 	 \proof The cancelations follow from the equalities
 	 \begin{align*}
 	 \sum_{1\leq i\leq 4} \left(\dellt (\{0,i\})-\dellt (\{0,i+1\})  \right)&=\dellt (\{0,1\})-\dellt (\{0,5\}) 
 	 \\
 	\sum_{1\leq i\leq 4} \left(\dellt (\{i,0\})-\dellt (\{i+1,0\})  \right)&=\dellt (\{1,0\})-\dellt (\{5,0\}) 
 	 \end{align*}
 	 and from the following one which uses the rules \eqref{rulepri}
 	 $$
 	 \sum_{1\leq i\leq 4}\left(\Delta (\{i,i+1\})-\Delta (\{i+1,i\})  \right)=0
 	 $$
 	 since \eqref{rulepri} shows that the following terms all vanish 
 	 $$
 	 \Delta(\{1,2\})- 	\Delta(\{4,3\}),  \ \Delta(\{2,3\})-\Delta(\{5,4\}),\ \Delta(\{3,4\})-\Delta(\{2,1\}), \ \Delta(\{4,5\})-\Delta(\{3,2\})
 	 $$
 	 We thus get the required formulas. \endproof 
 	 
 	 \begin{figure}[t!]
   \begin{minipage}{0.48\textwidth}
     \centering
      \includegraphics[width=.9\linewidth]{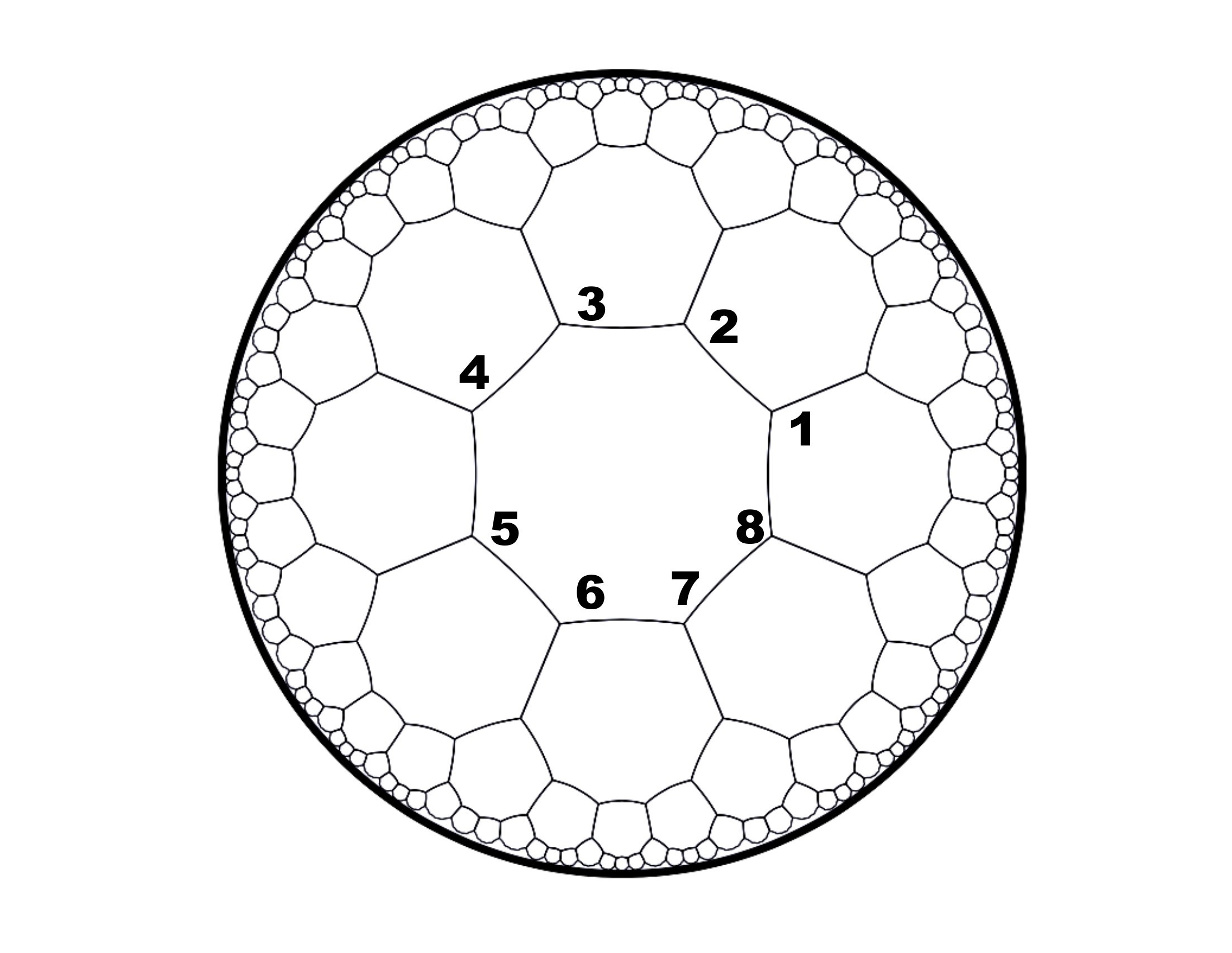}
          \caption{Surface of genus $2$.}\label{genus2}
   \end{minipage}\hfill
   \begin{minipage}{0.48\textwidth}
     \centering
    \includegraphics[width=.8\linewidth]{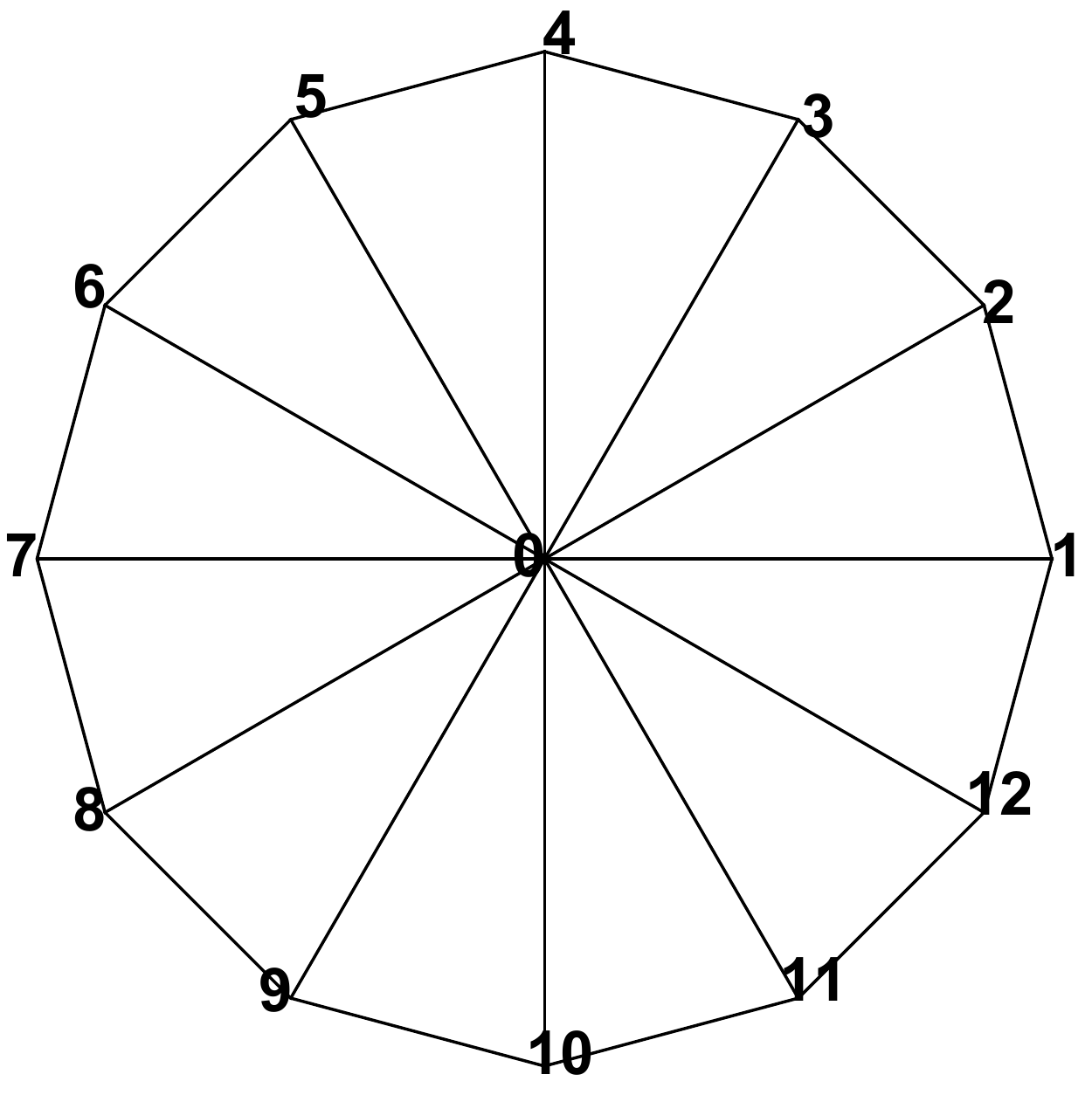}
     \caption{Domain for surface of genus $3$.}\label{octoexample}
   \end{minipage}
\end{figure}

 	 \subsection{Moore normalization for a Riemann surface of genus $g>1$}\label{block1}
 	 Let $g>1$ and $P$ be obtained (see Figure \ref{octoexample}) as the union of $g$ copies $P(w)=P_{(1+4w,5+4w)}$ for $0\leq w<g$, of the basic polygon  
 ${\rm Conv}(0,1,2,3,4,5)$ of Figure \ref{gromovnorm1}, where the side $(0, 5+4w)$ is common to $P(w)$ and $P(w+1)$ for $w<g-1$ and is common to $P(g-1)$ and $P(0)$ for $w=g-1$, while the external sides are identified pairwise as in $P$. The quotient of $P$ by these identifications is a surface $\Sigma(g)$ of genus $g$.
 	 	 \begin{lem}\label{chaingood3} The singular chain $$c=c_{(1,5)}+c_{(5,9)}+\ldots +c_{(1+4w,5+4w)}+\ldots +c_{(1+4(g-1),1)}$$ is closed and normalized, \ie one has $\partial_j c=0$ for $j\in \{0,1,2\}$. 	
 	\end{lem}
\proof By Lemma \ref{chaingood2} one gets for $0\leq w\leq g-1$, and with $5+4(g-1)\sim 1$,
 $$
 \partial_0 c_{(1+4w,5+4w)}=2 \dellt (\{0,1+4w\})-2 \dellt (\{0,5+4w\})-\dellt (\{1+4w,0\})+\dellt (\{5+4w,0\})
 $$
 which gives
 $$
 \partial_0 c=\partial_0 \sum_{0\leq w\leq g-1}  c_{(1+4w,5+4w)}=\sum_{0\leq w\leq g-1} \partial_0 c_{(1+4w,5+4w)}=0.
 $$	 
 	The same reasoning applies to show that $\partial_j c=0$ for $j\in \{1,2\}$. \endproof


 	\begin{lem}\label{chaingood4} The singular chain $c$ of Lemma \ref{chaingood3} represents the singular homology class $8[\Sigma]$ ($[\Sigma]=$ fundamental class of $\Sigma$). 	\end{lem} 
 	\proof The result follows since each chain $c(0,i,i+1)$ as in  \eqref{chaingood} is homologous to $8\dellt (\{0,i,i+1\})$ while the $4g$ triangles $\dellt (\{0,i,i+1\})$ for $1\leq i\leq 4g$ give a triangulation of $\Sigma$.\endproof 
 	\begin{lem}\label{chaingood5} The $\ell^1$-norm of the singular chain $c$ of Lemma \ref{chaingood3} is   $\leq 32 g$. 	\end{lem} 
 	 \proof This follows from the triangle inequality and the definition \eqref{chaingood} of the chain $c(0,i,j)$ whose $\ell^1$-norm is $\leq 8$. \endproof 
 	 \begin{thm}\label{chainthm} Let $\Sigma$ be a compact Riemann surface and $[\Sigma]$ its fundamental class in homology. Then $[\Sigma]$ belongs to the range of the canonical map $H_2(\Sigma,\Vert H\R \Vert_\lambda)\to H_2(\Sigma,\R)$ if and only if $\lambda$ is larger than the Gromov norm of $[\Sigma]$. 	 	
 	 \end{thm}
\proof The result  follows from Theorem \ref{propcompare1} if one shows that  the fundamental class $[\Sigma]$ fulfills the equality 
$$
\Vert [\Sigma]\Vert^{\rm nor}=\Vert [\Sigma]\Vert_1
$$
The inequality $\geq$ follows from \eqref{normN1}. Moreover, as recalled in section \ref{simplicialvol}, for a surface of genus $g$ the Gromov norm $\Vert [\Sigma]\Vert_1$ is equal to $4(g-1)$. Thus it remains to show that 
$\Vert [\Sigma]\Vert^{\rm nor}\leq 4(g-1)$. By applying Lemmas \ref{chaingood4} and \ref{chaingood5} one obtains the inequality $\Vert [\Sigma]\Vert^{\rm nor}\leq 4g$. One then applies a standard technique which is to use the same inequality for the covering space $\Sigma'$ of $\Sigma$ associated to an infinite cyclic subgroup of the fundamental group $\pi_1(\Sigma)$. The genus of a cyclic cover $\Sigma'$ of degree $n$ is $g'=n(g-1)+1$ since the Euler characteristic is multiplied by $n$. Thus the inequality $\Vert [\Sigma']\Vert^{\rm nor}\leq 4g'$ entails 
$$
n\Vert [\Sigma]\Vert^{\rm nor}\leq 4g'=4(n(g-1)+1).
$$
By passing to the limit when $n\to \infty$ one obtains the desired inequality $\Vert [\Sigma]\Vert^{\rm nor}\leq 4(g-1)$.\endproof

\end{document}